\documentclass[twoside,a4paper,reqno]{amsart}
\usepackage{graphicx}
\usepackage{amsmath}
\usepackage{amsfonts} % for \mathbb
\usepackage{amsthm}
\usepackage{amscd}
\usepackage{amsopn}
\usepackage{amstext}
\usepackage{amsxtra}
\usepackage{stmaryrd}
\usepackage{amssymb}
\usepackage{textcase} % For correct case of swedish chars
\usepackage[algo2e,ruled,linesnumbered]{algorithm2e} % for typesetting algorithms
\usepackage{bm} % for \boldsymbol
\numberwithin{equation}{section}
\numberwithin{figure}{section}
\theoremstyle{plain}
\newtheorem{thm}{Theorem}[section]
\newtheorem{lem}[thm]{Lemma}
\newtheorem{cl}[thm]{Corollary}
\newtheorem{prop}[thm]{Proposition}
\newtheorem{conj}[thm]{Conjecture}

\theoremstyle{definition}

\theoremstyle{remark}

\newcounter{algorithm2e}
\providecommand{\abs}[1]{\left\lvert #1 \right\rvert}

\providecommand{\ceil}[1]{\left\lceil #1 \right\rceil}

\providecommand{\set}[1]{\left\lbrace #1 \right\rbrace}
\providecommand{\gen}[1]{\left\langle #1 \right\rangle}

\renewcommand{\Pr}[1]{\operatorname{Pr}[ #1 ]}

\newcommand{\field}[1]{\mathbb{#1}}

\newcommand{\N}{\field{N}}
\newcommand{\Z}{\field{Z}}

\newcommand{\F}{\field{F}}
\newcommand{\PS}{\field{P}}

\newcommand{\MAGMA}{\textsc{Magma}}

\newcommand{\OV}{\mathcal{O}}

\DeclareMathOperator{\GL}{GL}

\DeclareMathOperator{\SL}{SL}

\DeclareMathOperator{\Sz}{Sz}
\DeclareMathOperator{\Sp}{Sp}

\DeclareMathOperator{\diag}{diag}
\DeclareMathOperator{\Tr}{Tr}
\DeclareMathOperator{\SLP}{\mathtt{SLP}}
\DeclareMathOperator{\Norm}{N}
\DeclareMathOperator{\Cent}{C}
\DeclareMathOperator{\Zent}{Z}

\newcommand{\OR}[1]{\operatorname{O} ( #1 )}

\title{Recognising the Suzuki groups in their natural representations}
\author{Henrik B\"a\"arnhielm}
\address{School of Mathematical Sciences \\ Queen Mary, University of London \\ Mile End Road \\ London E1 4NS \\ United Kingdom}
\urladdr{http://www.maths.qmul.ac.uk/\textasciitilde hb/}
\email{h.baarnhielm@qmul.ac.uk}

\begin{document}

%\begin{titlepage}
\begin{abstract}
  Under the assumption of a certain conjecture, for which there exists
  strong experimental evidence, we produce an efficient
  algorithm for constructive membership testing in the Suzuki groups
  $\Sz(q)$, where $q = 2^{2m + 1}$ for some $m > 0$, in their natural
  representations of degree $4$. It is a Las Vegas algorithm with
  running time $\OR{\log(q)}$ field operations, and a preprocessing
  step with running time $\OR{\log(q) \log{\log(q)}}$ field
  operations. The latter step needs an oracle for the discrete logarithm
  problem in $\F_q$.

  We also produce a recognition algorithm for $\Sz(q) = \gen{X}$. This is a Las
  Vegas algorithm with running time $\OR{\abs{X}^2}$
  field operations.

  Finally, we give a Las Vegas algorithm that, given $\gen{X}^h = \Sz(q)$
  for some $h \in \GL(4, q)$, finds some $g$ such that $\gen{X}^g =
  \Sz(q)$. The running time is $\OR{\log(q) \log{\log(q)} + \abs{X}}$
  field operations.

  Implementations of the algorithms are available for the computer
  system $\MAGMA$.
\end{abstract}

\maketitle
%\end{frontmatter}

\section{Introduction}

A goal of the \emph{matrix recognition project} is to develop
efficient algorithms for the study of subgroups of $\GL(d, q)$. The
classification due to Aschbacher (see \cite{aschbacher84}) provides
one framework for this, and the first aim is to develop an algorithm that finds
a composition series of a matrix group given by a set of generators.
It is possible to do this with a recursive
algorithm, and the recursion is described in \cite{crlg01}. However,
we still have to deal with the base cases, which are the finite simple
groups.

For each base case we need to perform parts of \emph{constructive
  recognition}. The simple group is given as $G = \gen{X}$ where $X \subseteq
\GL(d, q)$ for some $d, q$ and constructive recognition encompasses
the following problems:
\begin{enumerate}
\item The problem of \emph{recognition} or \emph{naming} of $G$, \emph{i.e.} decide the name of $G$, as in the classification of the finite simple groups.
\item The \emph{constructive membership} problem. Given $g \in
  \GL(d, q)$, decide whether or not $g \in G$, and if so express
  $g$ as a word (or $\SLP$, see Section \ref{section_slp}) in $X$.
\item Construct an isomorphism $\psi$ from $G$ to a \emph{standard copy} $H$ of $G$ such that $\psi(g)$ and $\psi^{-1}(h)$ can be computed efficiently for every $g \in G$ and $h \in H$. Sometimes this particular problem is what is meant by \lq \lq constructive recognition''.
\end{enumerate}

To find a composition series using \cite{crlg01}, we need only 
recognition and constructive membership, but the explicit isomorphisms
to a standard copy are also very useful. Given these, many problems,
including constructive membership, can be reduced to the standard copy.

This paper will consider the Suzuki groups $\Sz(q)$, $q = 2^{2m + 1}$
for $m > 0$, which is one of the infinite families of finite
simple groups. We will only consider the natural representation, which
has dimension $4$, and our standard copy will be $\Sz(q)$
defined in Section \ref{section:suzuki_theory}. 

In Section
\ref{section:constructive_membership} we solve the constructive
membership problem for $\Sz(q)$. In Section
\ref{section:standard_recognition} we solve the recognition problem
for $\Sz(q)$, \emph{i.e.} given $X \subseteq \GL(4, q)$ we give an
algorithm that decides whether or not $\gen{X} = \Sz(q)$. In Section
\ref{section:conjugation_problem} we consider these problems for conjugates of $\Sz(q)$.
Given $X \subseteq \GL(4, q)$ we give an algorithm that decides whether
or not $\gen{X}^h = \Sz(q)$ for some $h \in \GL(4, q)$. We also give
an algorithm that computes an isomorphism to $\Sz(q)$, by finding some
$g$ such that $\gen{X}^g = \Sz(q)$.

Other representations are dealt with in
\cite{sz_tensor_decompose}. The main objective of this paper is to
prove the following:

\begin{thm} \label{main_theorem} Assuming Conjecture
  \ref{conjecture_correctness}, and given a random element oracle for
  subgroups of $\GL(4, q)$ and an oracle for the discrete logarithm
  problem in $\F_q$, there exists a Las Vegas algorithm that, for each
  $X \subseteq \GL(4, q)$, with $q = 2^{2m + 1}$ for some $m > 0$,
  such that $\gen{X}^h = \Sz(q)$ for some $h \in \GL(4, q)$, finds $g
  \in \GL(4, q)$ such that $\gen{X}^g = \Sz(q)$ and solves the
  constructive membership problem for $\gen{X}$.  The algorithm has
  time complexity $\OR{\log(q)}$ field operations and also has a
  preprocessing step, which only needs to be executed once for a given
  $X$, with time complexity $\OR{\log(q) \log{\log(q)} + \abs{X}}$
  field operations. The discrete logarithm oracle is only needed in
  the preprocessing step.
\end{thm}
\begin{proof}
  Follows from Theorem \ref{thm_conj_problem}, Theorem
  \ref{thm_pre_step}, Theorem \ref{thm_element_to_slp} and Theorem
  \ref{thm_element_to_slp_complexity}.
\end{proof}

In Section \ref{section:implementation}, experimental evidence for Conjecture \ref{conjecture_correctness} is shown.

In constructive membership testing for $\Sz(q)$, the essential problem
is to find elements of even order. In this paper, this is achieved by
using the fact that $\Sz(q)$ acts doubly transitively on a certain set
$\OV \subseteq \PS^3(\F_q)$. After finding independent random elements in the
stabiliser of a point, which is done by finding elements that map one
point to another, it becomes easy to find elements of even order. This is
because the structure of the stabiliser of a point is known, and by
Proposition \ref{prop_frobenius} we can easily find elements of even
order in it.

For every cyclic subgroup $C$ of order $q - 1$, the
proportion of double cosets of $C$ in $\Sz(q)$ that contain an element
that maps one given point to another is high. The need to
consider double cosets rather than single cosets arises from the fact
that $\OV$ contains $q^2 + 1$ points, and most double
cosets have size $(q - 1)^2$. In the analogous problem for
$\SL(2, q)$ (see \cite{psl_recognition}), which acts on a set with $q
+ 1$ points, single cosets of a subgroup of order $q - 1$ are used.

One can view this as a process of applying permutation group
techniques on a set which is exponentially large in terms of the
input. Since $\OV$ has size $q^2 + 1$, we cannot explicitly write
down all its points and still have a polynomial time algorithm, and
therefore we cannot write down the elements of $\Sz(q)$ as
permutations. However, given two points
we can construct in polynomial time an element of $\Sz(q)$
that maps one point to the other, which is a typical permutation group
technique.

Implementations of the algorithms are available in $\MAGMA$ (see
\cite{magma}).

We are very grateful to the anonymous referee for the helpful advice and
the large number of comments. We also acknowledge John
Bray, Charles Leedham-Green, Eamonn O'Brien, Geoffrey Robinson, Maud
de Visscher and Robert Wilson for their help and encouragement.

\section{The simple Suzuki groups} \label{section:suzuki_theory}
We begin by defining our standard copy of the Suzuki group. Following \cite[Chapter $11$]{huppertIII}, let $\pi$ be the unique automorphism of
$\F_q$ such that $\pi^2(x) = x^2$ for every $x \in \F_q$, \emph{i.e.} $\pi(x) =
x^t$ where $t = 2^{m + 1}$. For $a, b \in \F_q$ and $c \in \F_q^{\times}$, define the following matrices.
\begin{equation}
S(a, b) = \begin{bmatrix}
1 & 0 & 0 & 0 \\
a & 1 & 0 & 0 \\
b & \pi(a) & 1 & 0 \\
a^2 \pi(a) + ab + \pi(b) & a \pi(a) + b & a & 1 
\end{bmatrix} 
\end{equation}
\begin{align}
M(c) &= \begin{bmatrix}
c^{1 + 2^m} & 0 & 0 & 0 \\
0 & c^{2^m} & 0 & 0 \\
0 & 0 & c^{-2^m} & 0 \\
0 & 0 & 0 & c^{-1 - 2^m}
\end{bmatrix} \\ 
T &= \begin{bmatrix}
0 & 0 & 0 & 1 \\
0 & 0 & 1 & 0 \\
0 & 1 & 0 & 0 \\
1 & 0 & 0 & 0 
\end{bmatrix}
\end{align}
By definition,
\begin{equation} \label{suzuki_def}
\Sz(q) = \gen{S(a, b), M(c), T \mid a, b \in  \F_q, c \in \F_q^{\times}}.
\end{equation}
If we define 
\begin{align}
\mathcal{F} &= \set{S(a, b) \mid a, b \in \F_q} \\
\mathcal{H} &= \set{M(c) \mid c \in \F_q^{\times}} 
\end{align}
then $\mathcal{F} \leqslant \Sz(q)$ with 
$\abs{\mathcal{F}} = q^2$ and $\mathcal{H} \cong \F_q^{\times}$ so
that $\mathcal{H}$ is cyclic of order $q - 1$. Moreover, we can write $M(c)$ as
\begin{equation}
M(c) = M^{\prime}(\lambda) = \begin{bmatrix}
\lambda^{t + 1} & 0 & 0 & 0 \\
0 & \lambda & 0 & 0 \\
0 & 0 & \lambda^{-1} & 0 \\
0 & 0 & 0 & \lambda^{-t-1}
\end{bmatrix}
\end{equation}
where $\lambda = c^{2^m}$.

The following result follows from \cite[Chapter $11$]{huppertIII}.
\begin{thm} \label{thm_suzuki_props}
\begin{enumerate}
\item The order of the Suzuki group is
\begin{equation} \label{suzuki_order}
\abs{\Sz(q)} = (q^2 + 1) q^2 (q - 1).
\end{equation}
\item For all $a, b, a^{\prime}, b^{\prime} \in \F_q$ and $\lambda \in \F_q^{\times}$ we have:
\begin{align}
S(a, b) S(a^{\prime}, b^{\prime}) &= S(a + a^{\prime}, b + b^{\prime} + a^t a^{\prime}) \label{matrix_id1} \\
S(a, b)^{M(\lambda)} &= S(\lambda a, \lambda^{t + 1} b). \label{matrix_id2}
\end{align}
\item There exists $\OV \subseteq \PS^3(\F_q)$ on which $\Sz(q)$ acts faithfully and doubly transitively, such that no nontrivial element of $\Sz(q)$ fixes
  more than $2$ points. This set is
\begin{equation} \label{ovoid_def}
\OV = \set{(1 : 0 : 0 : 0)} \cup \set{(a b + \pi(a) a^2 + \pi(b) : b : a : 1) \mid a,b \in \F_q}.
\end{equation}
\item The stabiliser of $P_{\infty} = (1 : 0 : 0 : 0) \in \OV$ is $\mathcal{F} \mathcal{H}$
and if $P_0 = (0 : 0 : 0 :1)$ then the stabiliser of $(P_{\infty},
P_0)$ is $\mathcal{H}$. 
\item $\Zent(\mathcal{F}) = \set{S(0, b) \mid b \in \F_q}$ and $\mathcal{FH}$ is a Frobenius group with Frobenius kernel $\mathcal{F}$.
\item The number of elements of order $q - 1$ is $\phi(q - 1) q^2 (q^2 + 1) / 2$, where $\phi$ is the Euler totient function.
\item Let $g \in G = \Sz(q)$. Then for every $x \in G$, $\Cent_G(g)
\cap \Cent_G(g)^x = \gen{1}$ if $\Cent_G(g) \neq \Cent_G(g)^x$.
\item $\Sz(q)$ has cyclic Hall subgroups $U_1$ and $U_2$ of orders $q \pm t + 1$.
\end{enumerate}
\end{thm}

From \cite[Chapter $11$, Remark $3.12$]{huppertIII} we also immediately obtain the following result.
\begin{thm} \label{sz_maximal_subgroups}
A maximal subgroup of $G = \Sz(q)$ is conjugate to one of the following subgroups.
\begin{enumerate}
\item The point stabiliser $\mathcal{F} \mathcal{H}$.
\item The normaliser $\Norm_G(\mathcal{H})$, which is dihedral of order $2(q - 1)$.
\item The normalisers $\mathcal{B}_i = \Norm_G(U_i)$ for $i = 1,2$. These satisfy $\mathcal{B}_i = \gen{U_i, t_i}$ where $u^{t_i} = u^q$ for every $u \in U_i$ and $[\mathcal{B}_i : U_i] = 4$.
\item $\Sz(s)$ where $q$ is a power of $s$.
\end{enumerate}
\end{thm}

If $G$ is a group acting on a set $\OV$ and $P \in \OV$, let $G_P \leqslant G$ denote the stabiliser of $P$ in $G$.

Let $\Sp(4, q)$ denote the standard copy of the symplectic group, preserving the following symplectic form:
\begin{equation} \label{standard_symplectic_form}
J = \begin{bmatrix}
0 & 0 & 0 & 1 \\
0 & 0 & 1 & 0 \\
0 & 1 & 0 & 0 \\
1 & 0 & 0 & 0
\end{bmatrix}.
\end{equation}

From \cite{ono62} and \cite[Chapter $3$]{raw04}, we know that the
elements of $\Sz(q)$ are precisely the fixed points of an automorphism
$\Psi$ of $\Sp(4, q)$; from \cite[Chapter $3$]{raw04}, computing
$\Psi(g)$ for some $g \in \Sp(4,q)$ amounts to taking a submatrix of
the exterior square of $g$ and then replacing each matrix entry $x$ by
$x^{2^m}$. Moreover, $\Psi$ is defined on $\Sp(4, F)$ for $F \geqslant \F_q$.

If $V$ is an $FG$-module for
some group $G$ and field $F$, with action $f : FG \times V \to
V$, and if $\Phi$ is an automorphism of $G$, denote by $V^{\Phi}$ the
$FG$-module which has the same elements as $V$ and where the action is given by $(g, v) \mapsto f(\Phi(g), v)$ for $g \in G$ and $v \in V^{\Phi}$, extended to $FG$ by linearity.

\begin{lem} \label{lem_steinberg_lang}
Let $G \leqslant \Sp(4, q)$ have natural module $V$ and assume that $V$ is absolutely irreducible. Then $G^h \leqslant \Sz(q)$ for some $h \in \GL(4, q)$ if and only if $V \cong V^{\Psi}$.
\end{lem}
\begin{proof}
  Assume $G^h \leqslant \Sz(q)$. Both $G$ and $\Sz(q)$ preserve the
  form \eqref{standard_symplectic_form}, and this form is unique up to a
  scalar multiple, since $V$ is absolutely irreducible. Therefore $h J
  h^T = \lambda J$ for some $\lambda \in \F_q^{\times}$. But if $\mu =
  \sqrt{\lambda^{-1}}$ then $(\mu h) J (\mu h)^T = J$, so that $\mu h \in
  \Sp(4, q)$. Moreover, $G^h = G^{\mu h}$, and hence we may assume
  that $h \in \Sp(4, q)$. Let $x = h \Psi(h^{-1})$ and observe that
  for each $g \in G$, $\Psi(g^h) = g^h$. It follows that
\begin{equation}
g^x = \Psi(h) g^h \Psi(h^{-1}) = \Psi(h g^h h^{-1}) = \Psi(g)
\end{equation}
so $V \cong V^{\Psi}$.

Conversely, assume that $V \cong V^{\Psi}$. Then there is some $h \in
\GL(4, q)$ such that for each $g \in G$ we have $g^h = \Psi(g)$. As above, since
both $G$ and $\Psi(G)$ preserve the form
\eqref{standard_symplectic_form}, we may assume that $h \in \Sp(4,
q)$. 

Let $K$ be the algebraic closure of $\F_q$. The Steinberg-Lang
Theorem (see \cite{steinberg_lang}) asserts that there exists $x \in
\Sp(4, K)$ such that $h = x^{-1} \Psi(x)$. It follows that
\begin{equation}
\Psi(g^{x^{-1}}) = \Psi(g)^{h^{-1} x^{-1}} = g^{x^{-1}}
\end{equation}
so that $G^{x^{-1}} \leqslant \Sz(q)$. Thus $G$ is conjugate in
$\GL(4, K)$ to a subgroup $S$ of $\Sz(q)$, and it follows from
\cite[Theorem $29.7$]{curtis_reiner}, that $G$ is conjugate to $S$ in
$\GL(4, q)$.
\end{proof}

\begin{lem} \label{lemma_double_coset}
If $H \leqslant G = \Sz (q)$ is a cyclic group of order $q - 1$ and $g \in G \setminus \Norm_G(H)$ then $\abs{HgH} = (q - 1)^2$.
\end{lem}
\begin{proof}
  Since $\abs{H} = q - 1$ it is enough to show that $H \cap H^g =
  \gen{1}$. By \cite[Chapter $11$]{huppertIII}, $H$ is conjugate to $\mathcal{H}$ and distinct conjugates of $\mathcal{H}$ intersect trivially.
\end{proof}

\begin{lem} \label{cl_totient_prop}
If $g \in G = \Sz(q)$ is uniformly random, then
\begin{equation}
\Pr{\abs{g} = q - 1} = \frac{\phi(q - 1)}{2(q - 1)} > \frac{1}{12\log{\log(q)}}
\end{equation}
and hence we expect to obtain an element of order $q - 1$ in
$\OR{\log\log{q}}$ random selections.
\end{lem}
\begin{proof}
The first equality follows immediately from Theorem \ref{thm_suzuki_props}. The
inequality follows from \cite[Section II.8]{totient_prop}.

Now let $\varepsilon = 1 / (12\log{\log(q)})$ and $\delta = \mathrm{e}^{-k}$
for some $k \in \N$. If we take uniformly random elements from $G$,
then the probability that we have not found an element of order $q -
1$ after $\ceil{\log{\delta} / \log{(1 - \varepsilon)}}$ consecutive tries is at most
$\delta$, and
\begin{equation}
\frac{\log{\delta}}{\log{(1 - \varepsilon)}} \approx \frac{k}{\varepsilon}
\end{equation}
which is $\OR{\log{\log(q)}}$, so the statement follows.
\end{proof}

\begin{lem} \label{lemma_fixing_elements}
The number of elements of $G = \Sz(q)$ that fix at least one point of $\OV$
is $q^2 (q - 1) (q^2 + q + 2) / 2$.
\end{lem}
\begin{proof}
  By \cite[Chapter $11$]{huppertIII}, if $g \in G$ fixes exactly one point,
  then $g$ is in a conjugate of $\mathcal{F}$ and if $g$ fixes two
  points then $g$ is in a conjugate of $\mathcal{H}$. This implies
  that there are $(\abs{\mathcal{F}} - 1)\abs{\OV}$ elements that fix
  exactly one point. Similarly, there are $\binom{\abs{\OV}}{2}(\abs{\mathcal{H}} - 1)$
  elements that fix exactly two points.
  
  Thus the number of elements that fix at least one point is
\begin{equation}
1 + (\abs{\mathcal{F}} - 1)\abs{\OV} + \binom{\abs{\OV}}{2}(\abs{\mathcal{H}} - 1) = \frac{q^2 ( q - 1) (q^2 + q + 2)}{2}.
\end{equation}
\end{proof}

\begin{lem} \label{lemma_conjugacy_classes}
Elements of odd order in $\Sz(q)$  that have the same trace are conjugate.
\end{lem}
\begin{proof}
  From \cite{suzuki62}, the number of conjugacy classes of
  non-identity elements of odd order is $q - 1$, and all elements of even
  order have trace $0$.  Observe that
\begin{equation}
S(0, b) T = \begin{bmatrix}
0 & 0 & 0 & 1 \\
0 & 0 & 1 & 0 \\
0 & 1 & 0 & b \\
1 & 0 & b & b^t 
\end{bmatrix}.
\end{equation}
Since $b$ can be any element of $\F_q$, so can $\Tr{(S(0, b) T)}$,
and this also implies that $S(0, b) T$ has odd order when $b \neq 0$.
Therefore there are $q - 1$ possible traces for
non-identity elements of odd order, and elements with different trace
must be non-conjugate, so all conjugacy classes must have different
traces.
\end{proof}

\section{Preliminaries}
We will now briefly discuss some general concepts that are needed later.

\subsection{Complexity} \label{complexity_general}
We shall be concerned with the time complexity of the algorithms
involved, where the basic operations are the field operations, and not
the bit operations. In our case, the matrix dimension will always be $4$, so all simple
arithmetic with matrices can be done using $\OR{1}$ field
operations, and raising a matrix to the $\OR{q}$ power can be done using $\OR{\log q}$ field operations using the standard method of
repeated squaring. We shall also assume an oracle for the
discrete logarithm problem for $\F_q$, so that this can be solved using $\OR{1}$ field operations. 

We will need to find an element of order $q - 1$. The order can be
computed using the algorithm of \cite{crlg95}. To obtain the
\emph{precise} order, this algorithm requires a factorisation of $q -
1$, otherwise it might return a multiple of the correct
order. However, it
suffices for our purposes to learn a
\emph{pseudo-order} of the element, which is a multiple of its
order, since it will suffice to find a nontrivial element of order dividing $q -
1$. Hence we avoid the requirement to factorise $q - 1$. The algorithm of \cite{crlg95} can also be used to obtain the pseudo-order, and for this it has time complexity $\OR{\log{(q)}\log{\log{(q)}}}$ field operations.

\subsection{Straight line programs} \label{section_slp}
For constructive membership testing, we want to express an element of a group $G = \gen{X}$ as
a word in $X$. Actually, it should be
a \emph{straight line program}, abbreviated to $\SLP$. If we express the elements as
words, the length of the words might be too large, requiring
exponential space complexity. 

An $\SLP$ is a data
structure for words, which ensures that subwords occurring multiple
times are computed only once. Formally, given a set of generators $X$,
an $\SLP$ is a sequence $(s_1, s_2, \dotsc, s_n)$ where
each $s_i$ represents one of the following
\begin{itemize}
\item an $x \in X$
\item a product $s_j s_k$, where  $j, k < i$
\item a power $s_j^n$ where $j < i$ and $n \in \Z$
\item a conjugate $s_j^{s_k}$ where $j, k < i$
\end{itemize}
so $s_i$ is either a pointer into $X$, a pair of pointers to earlier elements
of the sequence, or a pointer to an earlier element and an integer.

Thus to construct an $\SLP$ for a word, one starts by listing
pointers to the generators of $X$, and then builds up the word. To
evaluate the $\SLP$, go through the sequence and perform the
specified operations. Since we use pointers to the elements of $X$, we
can immediately evaluate the $\SLP$ on another set $Y$ of the same
size as $X$, by just changing the pointers so that they point to elements of $Y$.

\subsection{Random elements} 
Our analysis assumes that we can construct uniformly
distributed random elements of a group $G$ defined by a generating set
$X$. The polynomial time algorithm of \cite{babai91} produces nearly uniformly distributed random elements; an
alternative polynomial time algorithm is the \emph{product
replacement} algorithm of \cite{lg95}. We will assume that
we have a random element oracle, which produces a uniformly random element
using $\OR{1}$ field operations, and automatically gives it as an $\SLP$ in $X$. 

An important issue is the length of the $\SLP$s that are computed. The
length of the $\SLP$s must be polynomial, otherwise it would not be
polynomial time to evaluate them. We assume
that $\SLP$s of random elements have length $\OR{1}$.

\subsection{Las Vegas algorithms}
All the algorithms we consider are probabilistic of the type known as
\emph{Las Vegas algorithms}. This type of algorithm is discussed in
\cite[Section 25.8]{VonzurGathen03}, \cite[Section 1.3]{seress03} and
\cite[Section 3.2.1]{hcgt}. In short it is a probabilistic algorithm
with an input parameter $\varepsilon$ that either returns
\texttt{failure}, with probability at most $\varepsilon$, or otherwise
returns a correct result. The time complexity naturally depends on
$\varepsilon$.

We present Las Vegas algorithms as probabilistic algorithms that
either return a correct result, with probability bounded below by
$1/p(n)$ for some polynomial $p(n)$ in the size $n$ of the input, or
otherwise return \texttt{failure}.  By enclosing such an algorithm in
a loop that iterates $\ceil{\log \varepsilon / \log{(1 - 1/p(n))}}$
times, we obtain an algorithm that returns \texttt{failure} with
probability at most $\varepsilon$, and hence is a Las Vegas algorithm
in the above sense. Clearly if the enclosed algorithm is polynomial
time, the Las Vegas algorithm is polynomial time.

One can also enclose the algorithm in a loop that iterates until the
algorithm returns a correct result, thus obtaining a probabilistic
time complexity, and the expected number of iterations is then
$\OR{p(n)}$.

\section{Computing an element of a stabiliser}
\label{section:stabiliser_elements} As explained in the introduction,
in constructive membership testing for $\Sz(q)$ the essential problem
is to find an element of the stabiliser of a given point $P \in \OV$,
expressed as an $\SLP$ in our given generators $X$ of $G = \Sz(q)$.
The idea is to map $P$ to $Q \neq P$ by a random $g_1 \in G$, and
then compute $g_2 \in G$ such that $Pg_2 = Q$, so that $g_1
g_2^{-1} \in G_P$.

Thus the problem is to find an element that maps $P$ to $Q$, and the
idea is to look for it in double cosets of cyclic subgroups of order
$q - 1$. We first give an overview of the method.

Begin by selecting random $a, h \in G$ such that
$a$ has pseudo-order $q - 1$, and consider the equation
\begin{equation} \label{stab_eqn1}
P a^j h a^i = Q
\end{equation}
in the two indeterminates $i,j$. If we can solve this equation for $i$ and
$j$, thus obtaining positive integers $k, l$ such that $1 \leqslant k, l \leqslant q - 1$ and $P a^l h a^k = Q$, then we have an
element that maps $P$ to $Q$.

Since $a$ has order dividing $q - 1$, by \cite[Chapter
$11$]{huppertIII}, $a$ is conjugate to a matrix $M^{\prime}(\lambda)$ for some
$\lambda \in \F_q^{\times}$. This implies that we can diagonalise $a$
and obtain a matrix $x \in \GL(4, q)$ such that $M^{\prime}(\lambda)^x
= a$. It follows that if we define $P^{\prime} = P x^{-1}$,
$Q^{\prime} = Q x^{-1}$ and $g = h^{x^{-1}}$ then \eqref{stab_eqn1} is
equivalent to
\begin{equation} \label{stab_eqn2}
P^{\prime} M^{\prime}(\lambda)^j g M^{\prime}(\lambda)^i = Q^{\prime}.
\end{equation}

Now change indeterminates to $\alpha$ and $\beta$ by letting $\alpha = \lambda^j$ and $\beta = \lambda^i$, so that we obtain the following equation:
\begin{equation} \label{stab_eqn3}
P^{\prime} M^{\prime}(\alpha) g M^{\prime}(\beta) = Q^{\prime}.
\end{equation}
This determines four equations in $\alpha$ and $\beta$, and in Section
\ref{solve_magic_poly} we will describe how to find solutions for
them. A solution $(\gamma, \delta) \in \F_q^{\times} \times \F_q^{\times}$ determines 
$M^{\prime}(\gamma), M^{\prime}(\delta) \in \mathcal{H}$, and hence
also $c, d \in H = \mathcal{H}^x$. 

If $\abs{a} = q - 1$ then $\gen{a} = H$, so that there exists positive
integers $k$ and $l$ as above with $a^l = c$ and $a^k = d$, and these
integers can be found by computing discrete logarithms, since we also
have $\lambda^l = \gamma$ and $\lambda^k = \delta$. Hence we obtain a solution to \eqref{stab_eqn1} from the solution to \eqref{stab_eqn3}. If $\abs{a}$ is a
proper divisor of $q-1$, then it might happen that $c \notin \gen{a}$
or $d \notin \gen{a}$, but by Lemma \ref{cl_totient_prop} we know that
this is unlikely.

Thus the overall algorithm is as in Algorithm \ref{alg:stab_element}.
We show the time complexity of the algorithm in Section
\ref{section:trick_complexity} and prove that it is correct in Section
\ref{section:trick_correctness}.

\begin{algorithm2e} 
\dontprintsemicolon
\caption{\texttt{FindMappingElement}}
\SetKwFunction{RandomCyclicGroup}{RandomCyclicGroup}
\SetKwFunction{SolveEquation}{SolveEquation}
\SetKwFunction{Random}{Random}
\SetKwFunction{Diagonalise}{Diagonalise}
\SetKwFunction{Log}{DiscreteLog} 
\SetKw{aand}{and} 
\KwData{Generating set $X$ for $G = \Sz(q)$ and points $P \neq Q \in \OV$} 
\KwResult{An element $g$ of $G$, written as an $\SLP$ in $X$, such that $Pg = Q$}
\tcc{ Assumes the existence of a function \texttt{SolveEquation} that
  solves \eqref{stab_eqn3}, if possible. Also, assumes that the
  function \texttt{Random} returns an element as an $\SLP$ in $X$, and
  that \texttt{DiscreteLog} returns a positive integer if a discrete
  logarithm exists and $0$ otherwise.}
\Begin{
  $h := \Random(G)$ \;
  \tcc{ Find random element $a$ of pseudo-order $q - 1$}
  \BlankLine 
  $a := \Random(G)$ \;
  \If{$\abs{a} \mid q - 1$}
  {
    $(M^{\prime}(\lambda), x) := \Diagonalise(a)$ \;
    \tcc{ Now $M^{\prime}(\lambda)^x = a$} 
    \BlankLine
    \If{$\SolveEquation(h^{x^{-1}}, P x^{-1}, Q x^{-1})$} 
    { \label{alg:stab_element_solve_poly}
      Let $(\gamma, \delta)$ be a solution. \;
      $l := \Log(\lambda, \gamma)$ \;
      $k := \Log(\lambda, \delta)$ \;
      \If{$k > 0 \; \aand \; l > 0$ \label{alg:stab_element_discrete_log}}
      {
        \Return{$a^l  h  a^k$} \; \label{alg:stab_element_final_element} 
      }
    }
  }
  \Return{\texttt{fail}} 
}
\refstepcounter{algorithm2e}
\label{alg:stab_element}
\end{algorithm2e}

\subsection{Solving equation \eqref{stab_eqn3}} \label{solve_magic_poly} We
will now show how to obtain the solutions of \eqref{stab_eqn3}. It
might happen that there are no solutions, in which case the method
described here will detect this and return with failure.

By letting $P^{\prime} = (q_1 : q_2 :
q_3 : q_4)$, $Q^{\prime} = (r_1 : r_2 : r_3 : r_4)$ and $g = [g_{i,j}]$, we can write out \eqref{stab_eqn3}
and obtain
\begin{equation} \label{alpha_beta_eqns}
\begin{split}
(q_1 g_{1,1} \alpha^{t + 1} + q_2 g_{2,1} \alpha + q_3 g_{3,1} \alpha^{-1} + q_4 g_{4,1} \alpha^{-t-1}) \beta^{t + 1} &= C r_1 \\
(q_1 g_{1,2} \alpha^{t + 1} + q_2 g_{2,2} \alpha + q_3 g_{3,2} \alpha^{-1} + q_4 g_{4,2} \alpha^{-t-1}) \beta &= C r_2 \\
(q_1 g_{1,3} \alpha^{t + 1} + q_2 g_{2,3} \alpha + q_3 g_{3,3} \alpha^{-1} + q_4 g_{4,3} \alpha^{-t-1}) \beta^{-1} &= C r_3 \\
(q_1 g_{1,4} \alpha^{t + 1} + q_2 g_{2,4} \alpha + q_3 g_{3,4} \alpha^{-1} + q_4 g_{4,4} \alpha^{-t-1}) \beta^{-t - 1} &= C r_4
\end{split}
\end{equation}
for some constant $C \in \F_q$. Henceforth, we assume that $r_i \neq
0$ for $i = 1, \dotsc, 4$, since this is the difficult case, and also
extremely likely when $q$ is large, as can be seen from Proposition
\ref{prop_trick_general_case}. A method similar to the one described
in this section will solve \eqref{stab_eqn3} when some $r_i = 0$ and
Algorithm \ref{alg:stab_element} does not assume that all $r_i
\neq 0$.

\begin{prop} \label{prop_trick_general_case}
If $P^{\prime} = (p_1 : p_2 : p_3 : p_4) \in \OV^x$ is uniformly random, where $\OV^x = \set{Px \mid P \in \OV}$ for some $x \in \GL(4, q)$, then
\begin{equation}
\Pr{p_i \neq 0 \mid i = 1,\dotsc,4} \geqslant (1 - \frac{\sqrt{2q}}{q})^4.
\end{equation}
\end{prop}
\begin{proof}
Let $P^{\prime} = Px$ and $x = [x_{i,j}]$. If $P = (1 : 0 : 0 : 0)$ then $P^{\prime} = (x_{1,1} : x_{1,2} : x_{1,3} : x_{1,4})$ so clearly
\begin{multline}
\Pr{p_i = 0 \mid \text{some}\ i} \leqslant \frac{1}{\abs{\OV}} + (1 - \frac{1}{\abs{\OV}}) \\
(1 - \Pr{(a^{t + 2} + b^t + ab) x_{1,1} + x_{2,1}b + x_{3,1}a + x_{4,1} \neq 0 \mid a \neq 0, b \neq 0}^4).
\end{multline}
Now it follows that
\begin{multline}
\Pr{(a^{t + 2} + b^t + ab) x_{1,1} + x_{2,1}b + x_{3,1}a + x_{4,1} = 0 \mid a \neq 0, b \neq 0} = \\
= \sum_{k \in \F_q^{\times}} \Pr{(k^{t + 2} + b^t + kb) x_{1,1} + x_{2,1}b + x_{3,1}k + x_{4,1} = 0 \mid a = k, b \neq 0} \Pr{a = k} \leqslant \frac{t}{q}
\end{multline}
since in a field a polynomial of degree $t$ has at most $t$ roots. The result follows by observing that $t = \sqrt{2q}$.
\end{proof}

For convenience, we denote the expressions in the parentheses at
the left hand sides of \eqref{alpha_beta_eqns} as $K, L, M$ and $N$
respectively. Then if we let $C = L \beta r_2^{-1}$ we obtain three equations
\begin{equation} \label{alpha_beta_eqn1}
\begin{split} 
K \beta^t &= r_1 r_2^{-1} L \\
M \beta^{-2} &= r_3 r_2^{-1} L \\
N \beta^{-t-2} &= r_4 r_2^{-1} L 
\end{split}
\end{equation}
and in particular $\beta$ is a function of $\alpha$, since
\begin{equation} \label{beta_from_alpha}
\beta = \sqrt{L^{-1} M r_3^{-1} r_2}.
\end{equation}
By substituting the first two equations into the third in \eqref{alpha_beta_eqn1} we obtain
\begin{equation} \label{alpha_eqn1}
NK r_2 r_3 = r_1 r_4 ML
\end{equation}
and by raising the first equation to the $t$-th power and substituting into the second, we obtain
\begin{equation} \label{alpha_eqn2}
r_1 r_3^{t / 2} L^{1 + t/2} = r_2^{1 + t/2} M^{t/2} K.
\end{equation}
If instead we let $C = M \beta^{-1} r_3^{-1}$ and proceed similarly, we obtain two more equations
\begin{align}
N^t L r_3^{t + 1} &= M^{t + 1} r_2 r_4^t \label{alpha_eqn3} \\
N L^{t/2} r_3^{1 + t/2} &= M^{1 + t/2} r_4 r_2^{t/2}. \label{alpha_eqn4}
\end{align}
Now \eqref{alpha_eqn1}, \eqref{alpha_eqn2},
\eqref{alpha_eqn3} and \eqref{alpha_eqn4} are equations in $\alpha$
only, and by multiplying them by suitable powers of $\alpha$, they can
be turned into polynomial equations such that $\alpha$ only occurs to
the powers $ti$ for $i = 1, \dotsc, 4$ and to lower powers that are
independent of $t$. The suitable powers of $\alpha$ are $2t + 2$, $t +
t/2 + 2$, $2t + 3$ and $2t + t/2 + 2$, respectively.

Thus we obtain the following four equations.
\begin{equation} \label{alpha_poly_eqns}
\begin{split}
\alpha^{4t} c_1 + \alpha^{3t} c_2 + \alpha^{2t} c_3 + \alpha^t c_4 & = d_1 \\
\alpha^{4t} c_5 + \alpha^{3t} c_6 + \alpha^{2t} c_7 + \alpha^t c_8 & = d_2 \\
\alpha^{4t} c_9 + \alpha^{3t} c_{10} + \alpha^{2t} c_{11} + \alpha^t c_{12} & = d_3 \\
\alpha^{4t} c_{13} + \alpha^{3t} c_{14} + \alpha^{2t} c_{15} + \alpha^t c_{16} & = d_4 \\
\end{split}
\end{equation}
The $c_i$ and $d_j$ are polynomials in $\alpha$ with degree
independent of $t$, for $i = 1, \dotsc, 16$ and $j = 1, \dotsc, 4$
respectively, so \eqref{alpha_poly_eqns} can be considered a linear
system in the variables $\alpha^{nt}$ for $n = 1,\dotsc,4$, with
coefficients $c_i$ and $d_j$.  Now the aim is to obtain a single
polynomial in $\alpha$ of bounded degree. For this we need the
following conjecture.

\begin{conj} \label{conjecture_correctness} For every $P^{\prime} = P
  x^{-1}, Q^{\prime} = Q x^{-1}, g = h^{x^{-1}}$ where $P, Q \in \OV$,
  $h \in G$ and $x \in \GL(4, q)$, if we regard
  \eqref{alpha_poly_eqns} as simultaneous linear equations in the
  variables $\alpha^{nt}$ for $n = 1,\dotsc,4$, over the polynomial
  ring $\F_q[\alpha]$, then it has non-zero determinant.
\end{conj}

In other words, the determinant of the coefficients $c_i$ is not the
zero polynomial. We comment on the validity of Conjecture
\ref{conjecture_correctness} in Section \ref{section:implementation}.

\begin{lem} \label{lemma_magic_poly} Given $P^{\prime}, Q^{\prime}$ and $g$ as in Conjecture
  \ref{conjecture_correctness} and assuming Conjecture
  \ref{conjecture_correctness}, there exists a univariate polynomial
  $f(\alpha) \in \F_q[\alpha]$ of degree at most $60$, such that for every $(\gamma,
  \delta) \in \F_q^{\times} \times \F_q^{\times}$ that is a solution for $(\alpha, \beta)$ in \eqref{stab_eqn3} we have $f(\gamma) = 0$.
\end{lem}
\begin{proof}
  So far in this section we have shown that if we can solve \eqref{alpha_poly_eqns} we can also solve \eqref{stab_eqn3}. From the four equations of \eqref{alpha_poly_eqns} we can
  eliminate $\alpha^t$. We can solve for $\alpha^{4t}$ from the fourth
  equation, and substitute into the third,
  thus obtaining a rational expression with no occurrence of
  $\alpha^{4t}$. Continuing this way and substituting into the other
  equations, we obtain an expression for $\alpha^t$ in terms of the
  $c_i$ and the $d_i$ only. This can be substituted into any of the
  equations of \eqref{alpha_poly_eqns}, where $\alpha^{nt}$ for $n =
  1, \dotsc, 4$ is obtained by powering up the expression for
  $\alpha^t$. Thus we obtain a rational expression $f_1(\alpha)$ of degree
  independent of $t$. We now take $f(\alpha)$ to be the numerator of
  $f_1$.

  In other words, we think of the $\alpha^{nt}$ as independent
  variables and of \eqref{alpha_poly_eqns} as a linear system over
  these variables, with coefficients in $\F_q[\alpha]$. By Conjecture
  \ref{conjecture_correctness} we can solve this linear system.
  
  Two possible problems can occur: $f$ is identically zero or some of
  the denominators of the expressions for $\alpha^{nt}$, $n = 1,
  \dotsc, 4$ turn out to be $0$. However, Conjecture
  \ref{conjecture_correctness} rules out these possibilities. By
  Cramer's rule, the expression for $\alpha^t$ is a rational
  expression where the numerator is a determinant, so it consists of
  sums of products of $c_i$ and $d_j$. Each product consists of three
  $c_i$ and one $d_j$. By considering the calculations leading up to
  \eqref{alpha_poly_eqns}, it is clear that each of the products has
  degree at most $15$. Therefore the expression for $\alpha^{4t}$ and
  hence also $f(\alpha)$ has degree at most $60$.

We have only done elementary algebra to obtain $f(\alpha)$
from \eqref{alpha_poly_eqns}, and it is clear that
\eqref{alpha_poly_eqns} was obtained from \eqref{alpha_beta_eqns} by
elementary means only. Hence all solutions $(\gamma, \delta)$ to
\eqref{alpha_beta_eqns} must also satisfy $f(\gamma) = 0$, although
there may not be any such solutions, and $f(\alpha)$ may also have
other zeros.
\end{proof}

\begin{cl} \label{lemma_trick_time} Assuming Conjecture
  \ref{conjecture_correctness}, there exists a Las Vegas algorithm
  that, given $P^{\prime}, Q^{\prime}$ and $g$ as in Conjecture
  \ref{conjecture_correctness}, finds all $(\gamma, \delta) \in
  \F_q^{\times} \times \F_q^{\times}$ that are solutions of
  \eqref{stab_eqn3}. The algorithm has time complexity $\OR{\log{q}}$
  field operations.
\end{cl}
\begin{proof}
Let $f(\alpha)$ be the polynomial constructed in Lemma
\ref{lemma_magic_poly}.  To find all solutions to \eqref{stab_eqn3},
we find the zeros $\gamma$ of $f(\alpha)$, compute the
corresponding $\delta$ for each zero $\gamma$ using
\eqref{beta_from_alpha}, and check which pairs $(\gamma, \delta)$
satisfy \eqref{alpha_beta_eqns}. These pairs must be all solutions of
\eqref{stab_eqn3}.

The only work needed is simple matrix arithmetic, finding the roots of
a polynomial of bounded degree over $\F_q$, and raising matrices to
the power $t$, where $t \in \OR{q}$. Hence the time complexity is
$\OR{\log{q}}$ field operations and the algorithm is Las Vegas since
by \cite[Corollary 14.16]{VonzurGathen03} the algorithm for finding
the roots of $f(\alpha)$ is Las Vegas with this time complexity.
\end{proof}

By following the procedure outlined in Lemma \ref{lemma_magic_poly}, it is straightforward to obtain an expression for $f(\alpha)$, where
the coefficients are expressions in the entries of $g$, $P^{\prime}$
and $Q^{\prime}$, but we will not display it here, since it would take up too much space.

\subsection{Complexity} \label{section:trick_complexity}
\begin{thm} \label{thm_trick_time}
Given an oracle for the discrete logarithm problem in $\F_q$ and a random element oracle for $G$, the time
complexity of Algorithm \ref{alg:stab_element} is $\OR{\log(q)
\log{\log{(q)}}}$ field operations.
\end{thm}
\begin{proof}
Diagonalising a matrix uses $\OR{\log q}$ field operations,
since it involves finding the eigenvalues, \emph{i.e.} finding the roots of
a polynomial of constant degree over $\F_q$, see \cite[Corollary 14.16]{VonzurGathen03}.

Computing the pseudo-order of a matrix uses $\OR{\log (q)
  \log{\log{(q)}}}$ field operations, if we use the algorithm
described in \cite{crlg95}. From Corollary \ref{lemma_trick_time}, it follows that line \ref{alg:stab_element_solve_poly} uses $\OR{\log q}$ field
operations.

Finally, line \ref{alg:stab_element_final_element} uses $\OR{\log q}$
field operations, since the exponents are $\OR{q}$. We conclude that Algorithm
\ref{alg:stab_element} uses $\OR{\log{(q)} \log{\log{(q)}}}$ field operations.
\end{proof}

\subsection{Correctness} \label{section:trick_correctness} There are
two issues when considering the correctness of Algorithm
\ref{alg:stab_element}. Using the notation in the algorithm, we have
to show that \eqref{stab_eqn3} has a solution with high probability,
and that the integers $k$ and $l$ are positive with high probability.

The algorithm in Corollary \ref{lemma_trick_time} tries to find an
element in the double coset $\mathcal{H} g \mathcal{H}$, where $g =
h^{x^{-1}}$, and we will see that this succeeds with high probability
when $g \notin \Norm_G(\mathcal{H})$, which is very likely.

If the element $a$ has order precisely $q - 1$, then from the
discussion at the beginning of Section
\ref{section:stabiliser_elements}, we know that the integers $k$ and
$l$ will be positive. By Lemma \ref{cl_totient_prop} we know that it is likely that $a$ has order precisely $q - 1$ rather than just a divisor
of $q - 1$.

Hence it follows that Algorithm \ref{alg:stab_element} has high probability of success. We formalise this argument in the following results.

\begin{lem} \label{lemma_trick_correctness}
Assume Conjecture \ref{conjecture_correctness}. Let $G = \Sz(q)$ and let $P \in \OV$ and $a, h \in G$ be given, such that $\abs{a} = q - 1$. Let $Q \in \OV$ be uniformly random. If $h \notin \Norm_G(\gen{a})$, then
\begin{equation}
\frac{(q - 1)^2}{(q^2 + 1) \deg{f}}  \leqslant \Pr{Q \in P\gen{a} h \gen{a}} \leqslant \frac{(q - 1)^2}{q^2 + 1}
\end{equation}
where $f(\alpha)$ is the polynomial constructed in Lemma \ref{lemma_magic_poly}. If instead $h \in \Norm_G(\gen{a})$ then
\begin{equation}
\Pr{Q \in P\gen{a} h \gen{a}} = \frac{(q - 1) (q^2 - 1) + 2}{(q^2 + 1)^2}.
\end{equation}
\end{lem}
\begin{proof}
If $h \notin \Norm_G(\gen{a})$ then by Lemma \ref{lemma_double_coset}, $\abs{\gen{a} h \gen{a}} = (q -
1)^2$, and hence $\abs{ P\gen{a} h
\gen{a}} \leqslant (q - 1)^2$. 

On the other hand, for every $Q \in \OV$ we have 
\begin{equation} 
\abs{\set{(k_1, k_2) \mid k_1, k_2 \in \gen{a}, \, P k_1 h k_2 = Q}} \leqslant \deg{f}
\end{equation}
since this is the equation we consider in Section
\ref{solve_magic_poly}, and from Lemma \ref{lemma_magic_poly} we know
that all solutions must be roots of $f$. Thus $\abs{ P\gen{a} h
  \gen{a}} \geqslant \abs{\gen{a} h \gen{a}} / \deg{f}$. Since $Q$ is
uniformly random from $\OV$, and $\abs{\OV} = q^2 + 1$, the
result follows.

If $h \in \Norm_G(\gen{a})$ then $\gen{a} h
\gen{a} = h\gen{a}$ and $\abs{P h\gen{a}} = \abs{\gen{a}}$ if $\gen{a}$
does not fix $Ph$. By \cite[Chapter $11$]{huppertIII}, the number
of cyclic subgroups of order $q - 1$ is $\binom{\abs{\OV}}{2}$ and
$\abs{\OV} - 1$ such subgroups fix $Ph$. Moreover, if $\gen{a}$ fixes $Ph$
then $Ph \gen{a} = \set{Ph}$. Thus
\begin{multline}
%\begin{split}
\Pr{Q \in P \gen{a} h \gen{a}} = \Pr{Q \in P h \gen{a}} \Pr{P ha \neq Ph} + \\
+ \Pr{Q = Ph} \Pr{Ph a = Ph} = \frac{\abs{P h\gen{a}}}{\abs{\OV}} \left(1 - \frac{\abs{\OV} - 1}{\binom{\abs{\OV}}{2}}\right) + \frac{1}{\abs{\OV}} \frac{\abs{\OV} - 1}{\binom{\abs{\OV}}{2}} 
%\end{split}
\end{multline}
and the result follows.
\end{proof}

\begin{thm} \label{trick_alg}
Assuming Conjecture \ref{conjecture_correctness} and given a random element oracle for $G$ and an oracle for the discrete logarithm problem in $\F_q$, Algorithm \ref{alg:stab_element} is a Las Vegas algorithm that with probability
$s$ returns an element mapping $P$ to $Q$, where 
\begin{equation}
s > \frac{1}{12 \log{\log(q)} \deg{f}} + \OR{1/q}
\end{equation}
\end{thm}
\begin{proof}
  We use the notation from the algorithm. Let $g = h^{x^{-1}}$, $H =
  \mathcal{H}^x$, $P^{\prime} = P x^{-1}$ and $Q^{\prime} = Q x^{-1}$.
  Corollary \ref{lemma_trick_time} implies that line
  \ref{alg:stab_element_solve_poly} will succeed if $Q^{\prime} \in
  P^{\prime} \mathcal{H} g \mathcal{H}$. If $\abs{a} = q - 1$, then $H
  = \gen{a}$, and the previous condition is equivalent to $Q \in P
  \gen{a} h \gen{a}$.
  
  Moreover, if $\abs{a} = q - 1$ then line
  \ref{alg:stab_element_discrete_log} will always succeed. It might of course succeed when $\abs{a}$ is a
  proper divisor of $q - 1$, so it follows that $s$ satisfies the
  following inequality.
\begin{equation}
\begin{split}
s &\geqslant \Pr{\abs{a} = q - 1} (\Pr{h \in \Norm_G(\gen{a})} \Pr{Q \in P \gen{a} h \gen{a} \mid h \in \Norm_G(\gen{a})} + \\
&+\Pr{h \notin \Norm_G(\gen{a})} \Pr{Q \in P \gen{a} h \gen{a} \mid h \notin \Norm_G(\gen{a})})
\end{split}
\end{equation}
Since $h$ is uniformly random, using Theorem \ref{sz_maximal_subgroups} we obtain
\begin{equation}
\Pr{h \in \Norm_G(\gen{a})} = \frac{2 (q - 1)}{\abs{G}} = \frac{2}{q^2 (q^2 + 1)}
\end{equation}
From Lemma \ref{cl_totient_prop} and Lemma \ref{lemma_trick_correctness} we obtain
\begin{equation}
\begin{split}
s &\geqslant \frac{\phi(q - 1)}{2 (q - 1)} \left[\frac{(q - 1)^2}{(q^2 + 1) \deg{f}} - \frac{2}{q^2 (q^2 + 1)}\frac{(q - 1)^2}{(q^2 + 1)} + \frac{2}{q^2 (q^2 + 1)} \frac{2 + (q - 1) (q^2 - 1)}{(q^2 + 1)^2} \right] = \\
%&= \frac{(q^8 - 2 q^7 + 3 q^6 - 4 q^5 + (3 - 2 \deg{f}) q^4 + (6 \deg{f} - 2) q^3 + (2 - 6 \deg{f}) q^2 + 2 \deg{f} q - 4 \deg (f)) \phi(q - 1)}{2 \deg (f) (q - 1) q^2 (q^2 + 1)^3} \approx \frac{\phi(q - 1)}{2 \deg{f} q}
& = \frac{\phi(q - 1)}{2 (q - 1) \deg{f}} + \OR{1/q}
\end{split}
\end{equation}
and the probability of success follows from Lemma \ref{cl_totient_prop}.

Clearly if a solution is returned, it is correct, so the algorithm is
Las Vegas.
\end{proof}

\begin{cl} \label{thm_stab_elt} Assuming Conjecture
  \ref{conjecture_correctness} and given a random element oracle for subgroups of $\GL(4, q)$ and an oracle for the discrete logarithm problem in $\F_q$, there exists a Las Vegas algorithm that, given $X \subseteq \GL(4, q)$ such that $G = \gen{X} = \Sz(q)$
  and $P \in \OV$, finds a uniformly random $g \in G_P$, expressed as an $\SLP$ in $X$. The algorithm has time complexity $\OR{\log(q) \log \log(q)}$ field
  operations. If $s$ is as in Theorem \ref{trick_alg}, the probability
  of success is
\begin{equation}
s (1 - \frac{1}{\abs{\OV}}) > \frac{1}{12 \log{\log(q)} \deg{f}} + \OR{1/q}.
\end{equation}
\end{cl}
\begin{proof}
We compute $g$ as follows.
\begin{enumerate}
\item Find random $x \in G$. Let $Q = Px$ and return with failure if $P = Q$.
\item Use Algorithm \ref{alg:stab_element} to find $y \in G$ such that $Qy = P$.
\item Now $g = x y \in G_P$.
\end{enumerate}
Clearly this is a Las Vegas algorithm with probability of success as
stated. Moreover, the dominating term in the complexity is the call to
Algorithm \ref{alg:stab_element}, with time complexity given by
Theorem \ref{thm_trick_time}.

The element $g$ will be expressed as an $\SLP$ in $X$, since $x$ is
random and elements from Algorithm \ref{alg:stab_element} are
expressed as $\SLP$s.

Each call to Algorithm \ref{alg:stab_element} uses independent random
elements, so the double cosets under consideration are uniformly
random and independent. Therefore the elements returned by Algorithm
\ref{alg:stab_element} must be uniformly random. This implies that $g$
is uniformly random.
\end{proof}

\section{Constructive membership testing} \label{section:constructive_membership}

We will now give an algorithm for constructive membership
testing in $\Sz(q)$. Given a set of generators $X$, such that $G = \gen{X}
= \Sz(q)$, and given $g \in G$, we want to express $g$
as an $\SLP$ in $X$. We need the following result.

\begin{prop} \label{prop_frobenius}
If $g_1, g_2 \in \mathcal{FH}$ are uniformly random, then
\begin{equation}
\Pr{\abs{[g_1, g_2]} = 4} = 1 - \frac{1}{q - 1}.
\end{equation}
\end{prop}
\begin{proof}
  Let $A = \mathcal{FH} / \Zent(\mathcal{F})$. By Theorem
  \ref{thm_suzuki_props}, $[g_1, g_2] \in \mathcal{F}$ and has order
  $4$ if and only if $[g_1, g_2] \notin \Zent(\mathcal{F})
  \triangleleft \mathcal{FH}$. It therefore suffices to find the
  proportion of pairs $k_1, k_2 \in A$ such that $[k_1, k_2] = 1$.

  If $k_1 = 1$ then $k_2$ can be any element of $A$, which contributes $q(q
  - 1)$ pairs.  If $1 \neq k_1 \in \mathcal{F} / \Zent(\mathcal{F})
  \cong \F_q$ then $\Cent_A(k_1) = \mathcal{F} / \Zent(\mathcal{F})$,
  so we again obtain $q(q - 1)$ pairs. Finally, if $k_1 \notin
  \mathcal{F} / \Zent(\mathcal{F})$ then $\abs{\Cent_A(k_1)} = q - 1$
  so we obtain $q(q - 2)(q - 1)$ pairs. Thus we obtain $q^2 (q - 1)$ pairs
  from a total of $\abs{A \times A} =q^2 (q - 1)^2$ pairs, and the
  result follows.
\end{proof}

The algorithm for constructive membership testing has a preprocessing
step and a main step. The preprocessing step consists of finding
\lq \lq standard generators'' for $O_2(G_{P_{\infty}}) = \mathcal{F}$ and
$O_2(G_{P_0})$. In the case of $O_2(G_{P_{\infty}})$ the standard generators are defined as matrices $\set{S(a_i, x_i)}_{i = 1}^{n} \cup \set{S(0, b_i)}_{i =
1}^{n}$ for some unspecified $x_i \in \F_q$, such that $\set{a_1,
\dotsc, a_n}$ and $\set{b_1, \dotsc, b_n}$ form vector space bases of
$\F_q$ over $\F_2$ (so $n = \log_2{q} = 2m + 1$).

For every $a,b \in \F_q$, every matrix $S(a, b) \in G_{P_{\infty}}$ can be reduced to the identity by multiplying it by some of the standard generators of $O_2(G_{P_{\infty}})$, and similarly for $G_{P_0}$. The standard generators are therefore used in the main step to perform row operations in $G_{P_{\infty}}$ and $G_{P_0}$.

\begin{thm} \label{thm_pre_step} Assuming Conjecture
  \ref{conjecture_correctness} and given a random element oracle for $G$ and an oracle for the discrete logarithm problem in $\F_q$, the preprocessing step is a Las Vegas
  algorithm that finds standard generators for $O_2(G_{P_{\infty}})$
  and $O_2(G_{P_0})$. The preprocessing step has time complexity $\OR{\log(q)
    \log\log(q)}$ field operations. The probability of success is at
  least
\begin{equation}
r^4 \frac{\phi(q - 1)^2 (q - 2)^2}{(q - 1)^4} > \frac{1}{2^{10} 3^6 (\log \log(q))^6 (\deg{f})^4} + \OR{1/q}
\end{equation}
where $r$ is the success probability of the algorithm described in Corollary \ref{thm_stab_elt}. 
\end{thm}
\begin{proof}
The preprocessing step is the following:
\begin{enumerate}
\item Find random $a_1, a_2 \in G_{P_{\infty}}$ and $b_1, b_2
  \in G_{P_0}$ using the algorithm described in Corollary 
  \ref{thm_stab_elt}. Let $c_1 = [a_1, a_2]$, $c_2 = [b_1, b_2]$. 
\item Determine if $\abs{c_1} = \abs{c_2}
  = 4$, if $\abs{a_1}$ or $\abs{a_2}$ divides $q - 1$ and if $\abs{b_1}$ or $\abs{b_2}$ divides $q - 1$. Return with failure if any of these turn out to be false.

\item Let $d_1 \in \set{a_1, a_2}$ where $\abs{d_1}$ divides $q
- 1$, and let $d_2 \in \set{b_1, b_2}$ where $\abs{d_2}$
divides $q - 1$. Let $Y_{\infty} = \set{c_1, d_1}$ and $Y_0 =
\set{c_2, d_2}$. Diagonalise $d_1$ and obtain $M^{\prime}(\lambda) \in
G$, where $\lambda \in \F_q^{\times}$. Determine if $\lambda$ 
lies in a proper subfield of $\F_q$, and if so return with failure. Do similarly for $d_2$. 
  
\item As standard generators for $O_2(G_{P_{\infty}})$ we now take 
\begin{equation} \label{standard_gens}
L = \bigcup_{i = 1}^{2m + 1} \set{c_1^{d_1^i}, (c_1^2)^{d_1^i}}
\end{equation}
and similarly we obtain $U$ for $O_2(G_{P_0})$. 
\end{enumerate}

It follows from \eqref{matrix_id1} and \eqref{matrix_id2} that
\eqref{standard_gens} provides the standard generators for
$G_{P_{\infty}}$. These are expressed as $\SLP$s in $X$, since this is
true for the elements returned from the algorithm described in Corollary
\ref{thm_stab_elt}. 

By Corollary \ref{thm_stab_elt}, the first step succeeds with
probability $r^4$, and the random elements selected are uniformly
distributed and independent. Since $G_{P_{\infty}} = \mathcal{FH}$,
the proportion of elements of order $q - 1$ in $G_{P_{\infty}}$ is
$\phi(q - 1) / (q - 1)$, and similarly for $G_{P_0}$. Hence by
Proposition \ref{prop_frobenius}, the second step succeeds with
probability at least $(\phi(q - 1)^2 (q -2)^2) / (q - 1)^4$. If $\abs{d_1} = \abs{d_2} = q - 1$, the third step will also
succeed, since $\lambda$ will not lie in a proper subfield. Hence
$O_2(G_{P_{\infty}}) < \gen{Y_{\infty}} \leqslant G_{P_{\infty}}$ and
$\gen{Y_{\infty}} = G_{P_{\infty}}$ precisely when $d_1$ has order
$q-1$, and similarly for $Y_0$.

By the remark preceding the theorem, $L$ determines two sets of
field elements $\set{a_1, \dotsc, a_{2m+1}}$ and $\set{b_1, \dotsc,
b_{2m+1}}$. In this case each $a_i = a \lambda^i$ and $b_i = b
\lambda^{i(t + 1)}$, for some fixed $a,b \in \F_q^{\times}$, where $\lambda$
is as in the algorithm. Since $\lambda$ does not lie in a proper subfield, these sets form vector space bases of $\F_q$ over $\F_2$.

It then follows from Lemma \ref{cl_totient_prop} and Corollary
\ref{thm_stab_elt} that the probability of success of the
preprocessing step is as stated. Therefore the preprocessing step is
a Las Vegas algorithm.

We only determine if $d_1$ and $d_2$ have order dividing $q - 1$ in
order to obtain a polynomial time algorithm. To determine if $\lambda$ lies
in a proper subfield it suffices to determine if $\abs{\lambda} \mid
2^n - 1$ where $n$ is a proper divisor of $2m + 1$. Hence the
dominating term in the complexity is the computation of random
elements in the stabiliser, in the first step. The time complexity is
therefore the same as for the algorithm described in Corollary
\ref{thm_stab_elt}.
\end{proof}

Now we
consider the algorithm that expresses $g$ as an $\SLP$ in
$X$. It is given formally as Algorithm \ref{alg:main_alg}.

\begin{algorithm2e} 
\dontprintsemicolon
\caption{\texttt{ElementToSLP}}
\KwData{Standard generators $L$ for $G_{P_{\infty}}$ and
  $U$ for $G_{P_0}$. Matrix $g \in \gen{X} = G$.}
\KwResult{A $\SLP$ for $g$ in $X$.}  
\SetKwFunction{Random}{Random}
\Begin{
  $r := \Random(G)$ \;
  \If{$gr$ has an eigenspace $Q \in \OV$ \label{main_alg_find_point}}{
    Find $z_1 \in G_{P_{\infty}}$ using $L$ such that $Qz_1 =
  P_0$. \label{main_alg_row_op1} \;
  \tcc{ Now $(gr)^{z_1} \in G_{P_0}$.}
  \BlankLine
    Find $z_2 \in G_{P_0}$ using $U$ such that $(gr)^{z_1} z_2 = M^{\prime}(\lambda)$ for some $\lambda \in
  \F_q^{\times}$. \label{main_alg_row_op2} \;
    \tcc{ Express diagonal matrix as $\SLP$}
    \BlankLine
    $x := \Tr(M^{\prime}(\lambda))$ \;
    Find $h = [S(0, (x^t)^{1 / 4}), S(0, 1)^T]$ using $L \cup U$. \label{main_alg_row_op3} \;
    \tcc{ Now $\Tr{h} = x$.}
    \BlankLine
    Let $P_1, P_2 \in \OV$ be the fixed points of $h$. \;
    Find $a \in G_{P_{\infty}}$ using $L$ such that $P_1 a = P_0$. \label{main_alg_row_op4} \;
    Find $b \in G_{P_0}$ using $U$ such that $(P_2 a)b = P_{\infty}$. \label{main_alg_row_op5} \;
    \tcc{ Now $h^{ab} \in G_{P_{\infty}} \cap G_{P_0} = \mathcal{H}$,
  so $h^{ab} \in \set{M^{\prime}(\lambda)^{\pm 1}}$.}
    \BlankLine
    \eIf{$h^{ab} = M^{\prime}(\lambda)$}{
      Let $W$ be an $\SLP$ for $(h^{ab} z_2^{-1})^{z_1^{-1}} r^{-1}$. \label{main_alg_get_slp1} \;
      \Return{$W$} 
    }{
      Let $W$ be an $\SLP$ for $((h^{ab})^{-1} z_2^{-1})^{z_1^{-1}} r^{-1}$. \label{main_alg_get_slp2} \;
      \Return{$W$} 
      }
  }
  \Return{\texttt{fail}}
  }
\refstepcounter{algorithm2e}
\label{alg:main_alg}
\end{algorithm2e}

\begin{thm} \label{thm_element_to_slp}
Given a random element oracle for $G$, Algorithm \ref{alg:main_alg} is a Las Vegas algorithm with probability of success $1/2 + \OR{1/q}$.
\end{thm}
\begin{proof}
First observe that since $r$ is randomly chosen we obtain it as an $\SLP$. On line \ref{main_alg_find_point}
we check if $gr$ fixes a point, and from Lemma
\ref{lemma_fixing_elements} we see that
\begin{equation}
\Pr{gr \text{ fixes a point}} = \frac{q^2 + q + 2}{2 (q^2 + 1)} \approx \frac{1}{2}
\end{equation}

The elements found at lines \ref{main_alg_row_op1} and
\ref{main_alg_row_op2} can be computed using row operations, so we can
obtain them as $\SLP$s.

The element $h$ found at line \ref{main_alg_row_op3} clearly has trace
$x$, and it can be computed using row operations, so we obtain
it as an $\SLP$. From Lemma
\ref{lemma_conjugacy_classes} we know that $h$ is conjugate to
$M^{\prime}(\lambda)$ and therefore must fix $2$ points of
$\OV$. Hence lines
\ref{main_alg_row_op4} and \ref{main_alg_row_op5} make
sense, and the elements found can again be computed using row
operations and therefore we obtain them as $\SLP$s. 

The only elements in $\mathcal{H}$ that are conjugate to $h$ are
$M^{\prime}(\lambda)^{\pm 1}$, so clearly $h^{ab}$ must be one of them.

Finally, the elements that make up $W$ were found as $\SLP$s, and it is clear that if we evaluate $W$ we obtain
$g$. Hence the algorithm is Las Vegas and the theorem follows.
\end{proof}

\subsection{Complexity}
\begin{thm} \label{thm_element_to_slp_complexity}
  Given a random element oracle for $G$, Algorithm
  \ref{alg:main_alg} has time complexity $\OR{\log{q}}$ field
  operations, space complexity $\OR{\log^2{q}}$ and the length of the
  returned $\SLP$ is $\OR{\log{q}}$.
\end{thm}
\begin{proof}
From \eqref{standard_gens} we see that the number of standard
generators is $\OR{\log{q}}$, and each matrix uses $\OR{\log{q}}$ space, so the space complexity of the algorithm is $\OR{\log^2{q}}$.

This also immediately implies that the row operations performed at lines
\ref{main_alg_row_op1}, \ref{main_alg_row_op2},
\ref{main_alg_row_op3}, \ref{main_alg_row_op4} and
\ref{main_alg_row_op5} use $\OR{\log{q}}$ field operations.

Finding the fixed points of $h$, and performing the check at line
\ref{main_alg_find_point} only amounts to considering eigenspaces,
which uses $\OR{\log{q}}$ field operations. Thus the time complexity of
the algorithm is $\OR{\log{q}}$ field operations.

The $\SLP$s
returned from Algorithm \ref{alg:stab_element} have length $\OR{1}$,
and \eqref{standard_gens} implies that each standard generator also has
length $\OR{1}$. Hence because of our row operations, $W$ will have
length $\OR{\log{q}}$.
\end{proof}

\section{Recognition} \label{section:standard_recognition}
We now discuss how to recognise $\Sz(q)$. We are given a set $X \subseteq \GL(4, q)$ and we
want to decide whether or not $\gen{X} = \Sz(q)$, the group defined in
\eqref{suzuki_def}.

To do this, it suffices to determine if $X \subseteq \Sz(q)$ and if
$X$ does not generate a proper subgroup, \emph{i.e.} if $X$ is not
contained in a maximal subgroup. To determine if $g \in X$ is in
$\Sz(q)$, first determine if $\det(g) = 1$, then determine if $g$
preserves the symplectic form of $\Sp(4, q)$ and finally determine if
$g$ is a fixed point of the automorphism $\Psi$ of $\Sp(4, q)$,
mentioned in Section \ref{section:suzuki_theory}.

The recognition algorithm relies on the following result.
\begin{lem} \label{lemma_standard_recognition}
Let $H = \gen{X} \leqslant \Sz(q) = G$, where $X = \set{x_1, \dotsc, x_n}$
and let $C = \set{[x_i, x_j] \mid 1 \leqslant i < j \leqslant n}$ and $M$ be the natural module of $H$. Then $H = G$ if and only if the
following hold:
\begin{enumerate}
\item $M$ is an absolutely irreducible $H$-module.
\item $H$ is not conjugate in $\GL(4, q)$ to a subgroup of $\GL(4, r)$, where $q$ is a proper power of $r$.
\item $C \neq \set{1}$ and for every $c \in C \setminus \set{1}$ there exists $x \in X$ such that $[c, c^x] \neq 1$.
\end{enumerate}
\end{lem}
\begin{proof}
  By Theorem \ref{sz_maximal_subgroups}, the maximal subgroups of $G$
  that do not satisfy the first two conditions  are
  $\Norm_G(\mathcal{H})$, $\mathcal{B}_1$ and $\mathcal{B}_2$. For
  each, the derived group is contained in the normalised
  cyclic group, so all these maximal subgroups are metabelian. If $H$
  is contained in one of them and $H$ is not abelian, then $C \neq
  \set{1}$, but $[c, c^x] = 1$ for every $c \in C$ and $x \in X$ since
  the second derived group of $H$ is trivial. Hence the last
  condition is not satisfied.

Conversely, assume that $H = G$. Then clearly, the first two
conditions are satisfied, and $C \neq \set{1}$. Assume that the last
condition is false, so for some $c \in C \setminus \set{1}$ we have that
$[c, c^x] = 1$ for every $x \in X$. This implies that $c^x \in \Cent_G(c) \cap
\Cent_G(c)^{x^{-1}}$, and it follows from Theorem \ref{thm_suzuki_props}
that $\Cent_G(c) = \Cent_G(c)^{x^{-1}}$. Thus $\Cent_G(c) = \Cent_G(c)^g$ for all $g
\in G$, so $\Cent_G(c) \triangleleft G$, but $G$ is simple and we have a
contradiction.
\end{proof}

\begin{thm} \label{thm_standard_recognition}
  There exists a Las Vegas algorithm that, given $X \subseteq \GL(4, q)$, decides whether or not $\gen{X} = \Sz(q)$. Its time complexity is $\OR{\abs{X}^2}$
  field operations.
\end{thm}
\begin{proof}
The algorithm proceeds as follows.
\begin{enumerate}
\item Determine if every $x \in X$ is in $\Sz(q)$, and return \texttt{false} if not.

\item Determine if $\gen{X}$ is absolutely irreducible and if it is
not conjugate in $\GL(4, q)$ to a subgroup of $\GL(4, r)$, where $q$
is a proper power of $r$. Return \texttt{false} if any of these turn out to be false.

\item Using the notation of Lemma \ref{lemma_standard_recognition}, try to find $c \in C$ such that $c \neq 1$. Return \texttt{false} if it cannot be found.

\item If such $c$ can be found, and if $[c, c^x] \neq 1$ for some $x \in X$, then return \texttt{true}, else return \texttt{false}.

\end{enumerate}

From the discussion at the beginning of this section, the first step
is easily done using $\OR{\abs{X}}$ field operations. The MeatAxe (see \cite{meataxe} and \cite{better_meataxe}) can be used to determine if the
natural module is absolutely irreducible; the algorithm of \cite{smallerfield} can be used to determine if $\gen{X}$ is conjugate in $\GL(4, q)$ to a subgroup of $\GL(4, r)$, where $q$ is a proper power of $r$. Both these algorithms have time complexity $\OR{\abs{X}}$ field operations.

  The rest of the algorithm is a straightforward application of the last
  condition in Lemma \ref{lemma_standard_recognition}, except that it
  is sufficient to use the condition for one nontrivial commutator
  $c$. By Lemma
  \ref{lemma_standard_recognition}, if $[c, c^x] \neq 1$ then $\gen{X}
  = \Sz(q)$; but if $[c, c^x] = 1$, then $C_{\gen{X}}(c) \triangleleft
  \gen{X}$ and we cannot have $\Sz(q)$.

  It follows immediately that the time complexity of the
  algorithm is $\OR{\abs{X}^2}$ field operations. Since the MeatAxe is
  Las Vegas, this algorithm is also Las Vegas.
\end{proof}

\section{The conjugation problem} \label{section:conjugation_problem}

Given a conjugate $G$ of $\Sz(q)$ we describe an algorithm to
construct an isomorphism from $G$ to $\Sz(q)$ by finding a conjugating element. As one component, we need another recognition algorithm for $G$, since the one described in
Section \ref{section:standard_recognition} only works for the standard
copy of $\Sz(q)$. In \cite{general_recognition}, a general recognition algorithm
is described which could be used, but we prefer the
very fast algorithm described below, which works for
this special case.

\subsection{Recognition} \label{section:conjugate_recognition}
We want to determine if a given group $G = \gen{X} \leqslant \GL(4, q)$ is
a conjugate of $\Sz(q)$, without finding a conjugating element.
We consider carefully the subgroups of $\Sp(4, q)$ and rule out all
except those isomorphic to $\Sz(q)$. This relies on the fact that, up
to Galois automorphisms, $\Sz(q)$ has only one equivalence class of
faithful representations in $\GL(4,q)$ (see \cite{steinberg63}), so if
we can show that $G \cong \Sz(q)$ then $G$ is a conjugate of $\Sz(q)$.

\begin{thm}
  There exists a Las Vegas algorithm that, given $X \subseteq \GL(4,
  q)$, decides whether or not $\gen{X}^h = \Sz(q)$ for some $h \in
  \GL(4, q)$. The algorithm has time complexity $\OR{\abs{X}^2}$ field
  operations.
\end{thm}
\begin{proof}
  Let $G = \gen{X}$. The algorithm proceeds as follows.

\begin{enumerate}
\item Determine if $G$ is absolutely irreducible, using the MeatAxe, and return \texttt{false} if not.

\item Determine if $G$ preserves a non-zero symplectic form $M$. If so we
  conclude that $G$ is a subgroup of a conjugate of $\Sp(4, q)$, and if not then return \texttt{false}. This is essentially isomorphism testing of modules, which is described in \cite{meataxe}. Since $G$ is
  absolutely irreducible, the form is unique up to a scalar multiple.

\item Conjugate $G$ so that it preserves the form $J$. This amounts to finding a symplectic basis, \emph{i.e.} finding an
invertible matrix $X$ such that $X J X^{T} = M$, which is easily done.
Then $G^X$ preserves the form $J$ and thus $G^X \leqslant \Sp(4, q)$ so
that we can apply $\Psi$.

\item Determine if $V \cong V^{\Psi}$, where $V$ is the natural module for $G$ and $\Psi$ is the automorphism from Lemma \ref{lem_steinberg_lang}. If so we conclude that $G$ is a subgroup of some conjugate of $\Sz(q)$, and if not then return \texttt{false}.

\item Determine if $G$ is a proper subgroup of $\Sz(q)$, \emph{i.e.} if it is contained in a maximal subgroup. This can
be done using Lemma \ref{lemma_standard_recognition}. If so, then return \texttt{false}, else return \texttt{true}.
\end{enumerate}

The algorithms for finding a preserved form and for module isomorphism
testing are Las Vegas, with the same time complexity as the MeatAxe
(see \cite{meataxe} and \cite{better_meataxe}), which is
$\OR{\abs{X}}$ field operations since $G$ has constant degree. Hence
we obtain a Las Vegas algorithm, with the same time complexity as the
algorithm from Theorem \ref{thm_standard_recognition}.
\end{proof}

\subsection{Finding a conjugating element}
\label{section:conjugating_element} Now we assume that we are given $G
\leqslant \GL(4, q)$ such that $G^h = \Sz(q)$ for some $h \in \GL(4,
q)$, and we turn to the problem of finding some $g \in \GL(4, q)$ such
that $G^g = \Sz(q)$, thus obtaining an isomorphism from any conjugate of
$\Sz(q)$ to the standard copy. 

\begin{lem} \label{lem_find_ovoid_point} Given a random element oracle for subgroups of $\GL(4, q)$, there exists a Las Vegas
  algorithm that, given $X \subseteq \GL(4, q)$ such that $\gen{X}^h =
  \Sz(q)$ for some $h \in \GL(4, q)$, finds a point $P \in
  \OV^{h^{-1}} = \set{Q h^{-1} \mid Q \in \OV}$. The algorithm has time complexity $\OR{\log q}$ field operations.
\end{lem}
\begin{proof}
  Clearly $\OV^{h^{-1}}$ is the set on which $\gen{X}$ acts doubly
  transitively. For a matrix $M^{\prime}(\lambda) \in \Sz(q)$ we see
  that the eigenspaces corresponding to the eigenvalues $\lambda^{\pm
    (t + 1)}$ will be in $\OV$. Moreover, every element of order
  dividing $q-1$ in every conjugate $G$ of $\Sz(q)$ will have
  eigenvalues of the form $\mu^{\pm (t + 1)}$, $\mu^{\pm 1}$ for some
  $\mu \in \F_q^{\times}$, and the eigenspaces corresponding to $\mu^{\pm (t +
    1)}$ will lie in the set on which $G$ acts doubly transitively.

  Hence to find a point $P \in \OV^{h^{-1}}$ it suffices to find a
  random $g \in \gen{X}$ of order dividing $q - 1$, which is easy by
  Lemma \ref{cl_totient_prop}, and then find the eigenspaces of $g$.

  Clearly this is a Las Vegas algorithm that uses $\OR{\log q}$ field operations.
\end{proof}

\begin{lem} \label{lem_diagonal_conj}
  There exists a Las Vegas algorithm that, given $X \subseteq \GL(4, q)$ such that $\gen{X}^d = \Sz(q)$ where
  $d = \diag(d_1, d_2, d_3, d_4) \in \GL(4, q)$, finds a diagonal matrix $e
  \in \GL(4, q)$ such that $\gen{X}^e = \Sz(q)$, using
  $\OR{\abs{X} + \log{q}}$ field operations.
\end{lem}
\begin{proof}
Let $G = \gen{X}$. Since $G^d = \Sz(q)$, $G$ must preserve the symplectic form
\begin{equation}
K = d J d = \begin{bmatrix}
0 & 0 & 0 & d_1 d_4 \\
0 & 0 & d_2 d_3 & 0 \\
0 & d_2 d_3 & 0 & 0 \\
d_1 d_4 & 0 & 0 & 0 
\end{bmatrix}
\end{equation}
where $J$ is given by \eqref{standard_symplectic_form}. Using \cite{meataxe}, we can find this form, which is determined up
to a scalar multiple. Hence the diagonal matrix $e = \diag(e_1, e_2, e_3, e_4)$
that we want to find is also determined up to a scalar multiple (and up to multiplication by a diagonal matrix in $\Sz(q)$).

Since $e$ must take $J$ to $K$, we must have
$K_{1, 4}= d_1 d_4 = e_1 e_4$ and $K_{2, 4} = d_2 d_3 = e_2 e_3$. The
matrix $e$ is determined up to a scalar multiple, so we can choose
$e_4 = 1$ and $e_1 = K_{1, 4}$. Hence it only remains to determine
$e_2$ and $e_3$.

To conjugate $G$ into $\Sz(q)$ we must have $Pe \in \OV$ for every point
$P \in \OV^{d^{-1}}$, which is the set on which $G$ acts doubly
transitively. By Lemma \ref{lem_find_ovoid_point}, we can find $P = (p_1 : p_2 : p_3 : 1) \in \OV^{d^{-1}}$, and the
condition $P e = (p_1 K_{1, 4} : p_2 e_2 : p_3 e_3 : 1) \in \OV$ is
given by \eqref{ovoid_def} and amounts to
\begin{equation} \label{conjugate_final_eq1}
p_2 p_3 K_{2, 3} + (p_2 e_2)^t + (p_3 e_3)^{t + 2} - p_1 K_{1, 4} = 0
\end{equation}
which is a polynomial equation in the two variables $e_2$ and
$e_3$. 

Notice that we can consider $e_2^{t}$ to be the
variable, instead of $e_2$, since if $x = e_2^{t}$, then $e_2 =
\sqrt{x^t}$. Similarly, we can let $e_3^{t + 2}$ be the variable
instead of $e_3$, since if $y = e_3^{t + 2}$ then $e_3 = y^{1 - t/2}$.
Thus instead of \eqref{conjugate_final_eq1} we obtain a linear equation
\begin{equation} \label{conjugate_final_eq2}
p_2^t x + p_3^{t + 2} y = p_1 K_{1, 4} - p_2 p_3 K_{2, 3}
\end{equation}
in the variables $x, y$. Thus the complete algorithm for finding $e$
proceeds as follows.
\begin{enumerate}
\item Find the form $K$ that is preserved by $G$, using \cite{meataxe}.
\item Find $P, Q \in \OV^{d^{-1}}$ using Lemma \ref{lem_find_ovoid_point}. 
\item Let $P = (p_1 : p_2 : p_3 : p_4)$ and $Q = (q_1 : q_2 : q_3 : q_4)$. Determine if the following linear system in the variables $x$ and $y$ is singular, and if so return with failure.
\begin{equation}
\begin{split}
p_2^t x + p_3^{t + 2} y &= p_1 K_{1, 4} - p_2 p_3 K_{2, 3} \\
q_2^t x + q_3^{t + 2} y &= q_1 K_{1, 4} - q_2 q_3 K_{2, 3}
\end{split}
\end{equation}
\item Let $(\alpha, \beta)$ be a solution to the linear system. The diagonal matrix $e = \diag(K_{1, 4}, \sqrt{\alpha^t}, \beta^{1 - t/2}, 1)$ now satisfies that $G^e = \Sz(q)$.
\end{enumerate}
By Lemma \ref{lem_find_ovoid_point} and \cite{meataxe}, this
is a Las Vegas algorithm that uses $\OR{\abs{X} + \log q}$ field
operations.
\end{proof}

\begin{lem} \label{lem_conjugate_to_digaonal}
  There exists a Las Vegas algorithm that, given subsets $X$, $Y_P$
  and $Y_Q$ of $\GL(4, q)$ such that $O_2(G_P) < \gen{Y_P} \leqslant
  G_P$ and $O_2(G_Q) < \gen{Y_Q} \leqslant G_Q$, respectively, where
  $\gen{X} = G$, $G^h = \Sz(q)$ for some $h \in \GL(4, q)$ and $P, Q
  \in \OV^{h^{-1}}$, finds $k \in \GL(4, q)$ such that $(G^k)^d =
  \Sz(q)$ for some diagonal matrix $d \in \GL(4, q)$.  The algorithm
  has time complexity $\OR{\abs{X}}$ field operations.
\end{lem}
\begin{proof}

  Notice that the natural module $V = \F_q^4$ of $\mathcal{F}
  \mathcal{H}$ is uniserial with four non-zero submodules, namely $V_i
  = \set{(v_1, v_2, v_3, v_4) \in \F_q^4 \mid v_j = 0, j > i}$ for $i
  = 1, \dotsc, 4$. Hence the same is true for $\gen{Y_P}$ and
  $\gen{Y_Q}$ (but the submodules will be different) since they lie in
  conjugates of $\mathcal{F} \mathcal{H}$.

  Now the algorithm proceeds as follows.
  
\begin{enumerate}
\item Let $V = \F_q^4$ be the natural module for $\gen{Y_P}$ and
  $\gen{Y_Q}$. Find composition series $V = V^P_4 \supset V^P_3
  \supset V^P_2 \supset V^P_1$ and $V = V^Q_4 \supset V^Q_3 \supset
  V^Q_2 \supset V^Q_1$ using the MeatAxe.

\item Let $U_1 = V_1^P$, $U_2 = V_3^P \cap V_2^Q$, $U_3 = V_2^P \cap V_3^Q$ and $U_4 = V_1^Q$. For each $i = 1, \dotsc, 4$, choose $u_i \in U_i$.

\item Now let $k$ be the matrix such that $k^{-1}$ has $u_i$ as row $i$, for $i = 1, \dotsc, 4$.
\end{enumerate}

We now motivate the second step of the algorithm. Let $(M)_i$ denote the $i$-th row of a matrix $M$, and let $V_i^P$ and $V_i^Q$ be as in the algorithm. 

  We may assume that $Y_P = \set{x, y}$, $Y_Q = \set{u, v}$
  where $\abs{x} = \abs{u} = 4$ and both $\abs{y}$ and $\abs{v}$
  divide $q - 1$ (and $y$ and $v$ are nontrivial). 

  There exists $g^{\prime}
  \in \Sz(q)$ such that $P h g^{\prime} = P_{\infty}$ and $Q h
  g^{\prime} = P_0$, since $\Sz(q)$ acts doubly transitively on $\OV$.
  If we let $z = hg^{\prime}$, then $\gen{Y_P}^z$ and $\gen{Y_Q}^z$ consist of
  lower and upper triangular matrices, respectively. Hence there
exist $a_1, b_1 \in \F_q$ such that $x = S(a_1, b_1)^{z^{-1}}$, and then $V_1^P =
\gen{(x)_1} = \gen{(S(a_1, b_1))_1 z^{-1}} = V_1$. But $(S(a_1, b_1))_1 z^{-1} = (z^{-1})_1$
so by choosing some non-zero vector in $V_1^P$ we obtain a scalar multiple
of the first row of $z^{-1}$. Similarly, there exist $a_2, b_2 \in \F_q$ such
that $u = (S(a_2, b_2)^{T})^{z^{-1}}$, and $V_1^Q = \gen{(u)_4} = \gen{(S(a_2,
  b_2)^{T})_4 z^{-1}}$, where $S(a_2, b_2)^{T}$ is the transpose of $S(a_2, b_2)$. But $(S(a_2, b_2)^{T})_4 z^{-1} = (z^{-1})_4$ so by choosing some
non-zero vector in $V_1^Q$ we obtain a scalar multiple of the fourth row
of $z^{-1}$.

Note that $\dim V^P_3 \cap
V^Q_2 = 1$ and $\dim V^P_2 \cap V^Q_3 = 1$, and by choosing non-zero
vectors from these we obtain scalar multiples of the second and third
rows of $z^{-1}$, respectively. 

Thus the matrix $k$ found in the algorithm satisfies that $z = kd$ for
some diagonal matrix $d \in \GL(4, q)$. Since $\Sz(q) = G^h = G^z =
(G^k)^d$, the algorithm returns a correct result, and it is Las Vegas
because the MeatAxe is Las Vegas (see \cite{meataxe} and \cite{better_meataxe}). Clearly the time complexity is the
same as the MeatAxe, so the algorithm uses $\OR{\abs{X}}$ field
operations.
\end{proof}

\begin{thm} \label{thm_conj_problem}
  Assuming Conjecture \ref{conjecture_correctness} and given a random element oracle for subgroups of $\GL(4, q)$, there exists a Las
  Vegas algorithm that, given $X \subseteq \GL(4, q)$ such that
  $\gen{X}^h = \Sz(q)$ for some $h \in \GL(4, q)$, finds $g \in \GL(4,
  q)$ such that $\gen{X}^g = \Sz(q)$. The algorithm has time complexity
  $\OR{\log(q) \log{\log(q)} + \abs{X}}$ field operations.
\end{thm}
\begin{proof}
  Let $G = \gen{X}$. First note that $g$ is determined up to
  multiplication by an element of $\Sz(q)$, so we will find $g$ such
  that $h g^{\prime} = g $ where $g^{\prime} \in \Sz(q)$. 

  The algorithm described in Corollary \ref{thm_stab_elt} works
  equally well for a conjugate of $\Sz(q)$, so we can find generators
  for a stabiliser of a point in $G$, using the algorithm described in
  Theorem \ref{thm_pre_step}. In this case we do not need
  the elements as $\SLP$s, so a discrete log oracle is not necessary.

\begin{enumerate}
\item Find points $P, Q \in \OV^{h^{-1}}$ using Lemma
  \ref{lem_find_ovoid_point}. Return with failure if $P = Q$.
\item Find generating sets $Y_P$ and $Y_Q$ such that $O_2(G_P) <
  \gen{Y_P} \leqslant G_P$ and $O_2(G_Q) < \gen{Y_Q} \leqslant G_Q$
  using the first three steps of the algorithm from the proof of Theorem
  \ref{thm_pre_step}.
\item Find $k \in \GL(4, q)$ such that $(G^k)^d = \Sz(q)$ for some diagonal matrix $d \in \GL(4, q)$, using Lemma \ref{lem_conjugate_to_digaonal}.
\item Find a diagonal matrix $e$ using Lemma \ref{lem_diagonal_conj}.
\item Now $g = ke$ satisfies that $G^g = \Sz(q)$.
\end{enumerate}

Be Lemma \ref{lem_find_ovoid_point}, \ref{lem_conjugate_to_digaonal} and \ref{lem_diagonal_conj}, and the proof of Theorem \ref{thm_pre_step}, this is a Las Vegas algorithm with time complexity as stated.
\end{proof}

\section{Implementation and performance} \label{section:implementation}

An implementation of the algorithms described here is available in
$\MAGMA$. The implementation uses the existing $\MAGMA$
implementations of the algorithms described in \cite{crlg95}, \cite{lg95},
\cite{smallerfield}, \cite{meataxe} and \cite[Corollary 14.16]{VonzurGathen03}.

\begin{figure}[ht]
\includegraphics[scale=0.65]{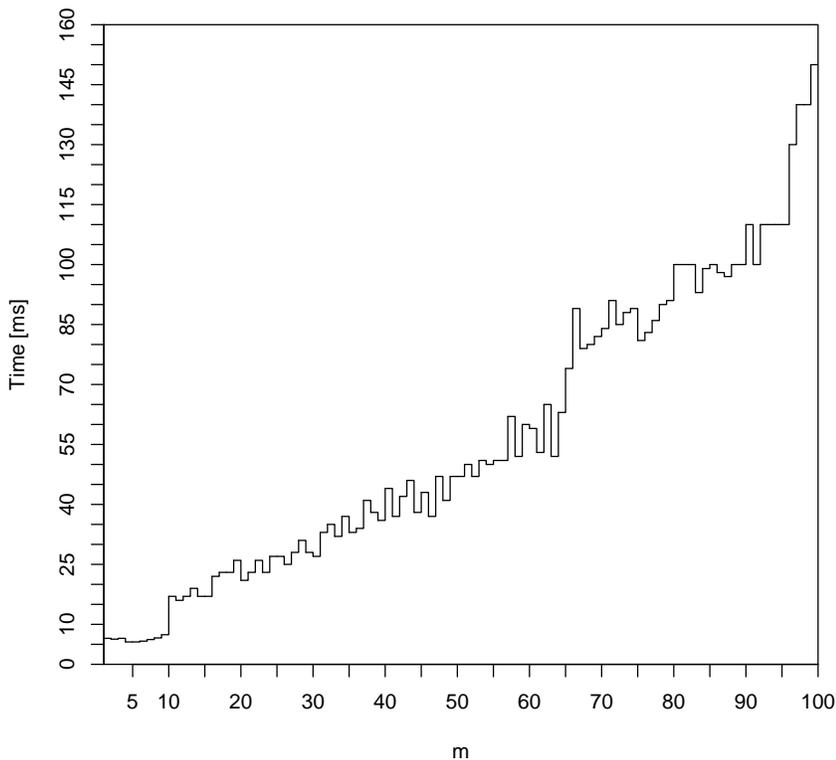}
\caption{Benchmark of recognition}
\label{fig:recognition_benchmark}
\end{figure}

A benchmark of the recognition algorithm described in Section
\ref{section:conjugate_recognition}, for various field sizes $q = 2^{2m
  + 1}$, is given in Figure \ref{fig:recognition_benchmark}. For each
field size, $200$ random conjugates of $\Sz(q)$ were recognised and the
average running time for each call is displayed.

\begin{figure}[ht]
\includegraphics[scale=0.65]{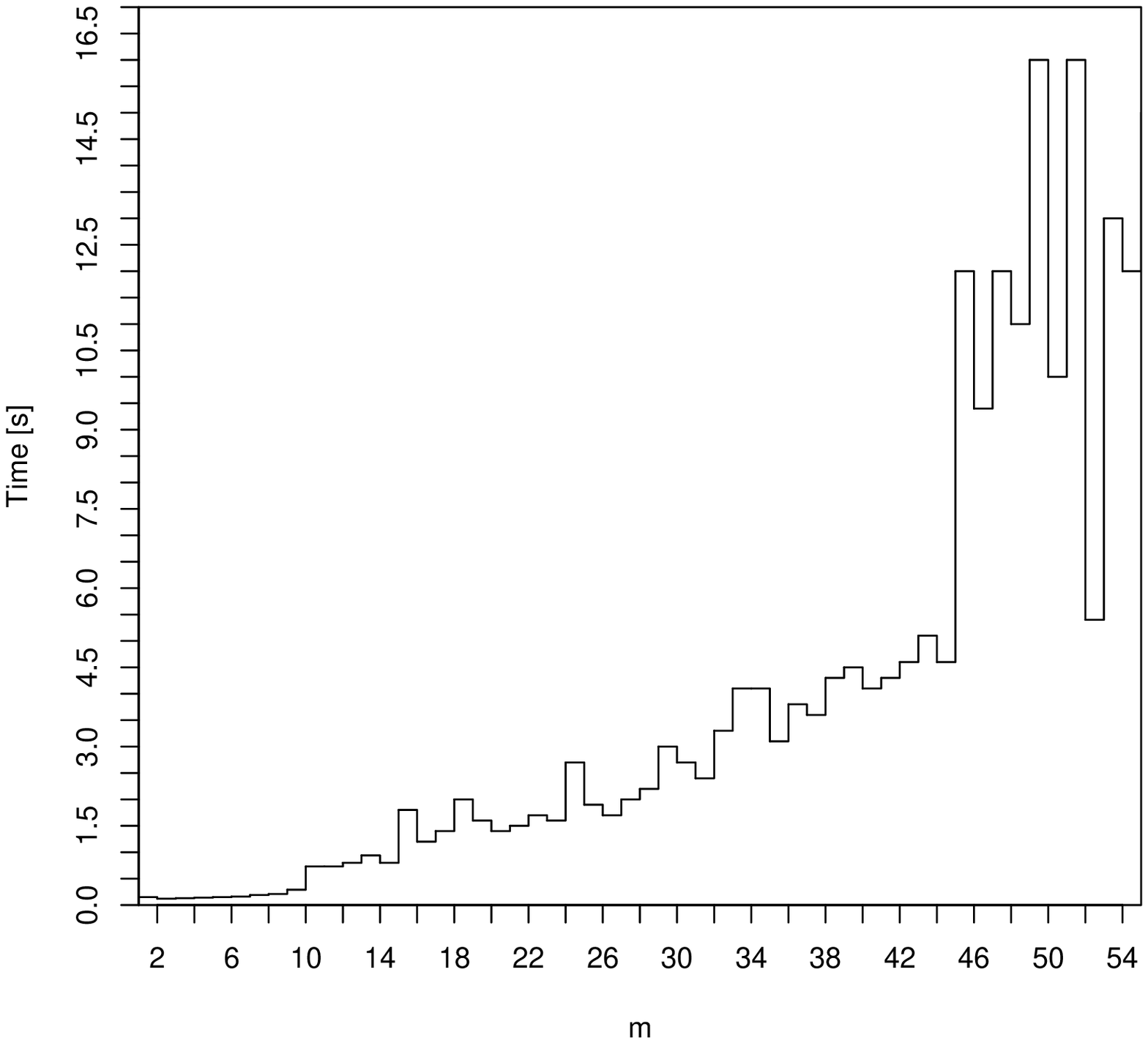}
\caption{Benchmark of conjugation}
\label{fig:conjugation_benchmark}
\end{figure}

A benchmark of the conjugation algorithm described in Section
\ref{section:conjugating_element}, for various field sizes $q = 2^{2m
+ 1}$, is given in Figure \ref{fig:conjugation_benchmark}. For each
field size, $100$ random conjugates of $\Sz(q)$ were considered and
a conjugating element found. The average running time for each call
is displayed.

The constructive membership and
conjugation algorithms both need to compute generating sets of
stabilisers, so they depend on Algorithm \ref{alg:stab_element}.
Therefore our implementation depends on the $\MAGMA$ implementation of
discrete log. Since we are in characteristic $2$, there
is a specialised algorithm for discrete log, \emph{Coppersmith's
  algorithm} (see \cite{cop84}), which is implemented in $\MAGMA$.

\begin{figure}[ht]
\includegraphics[scale=0.65]{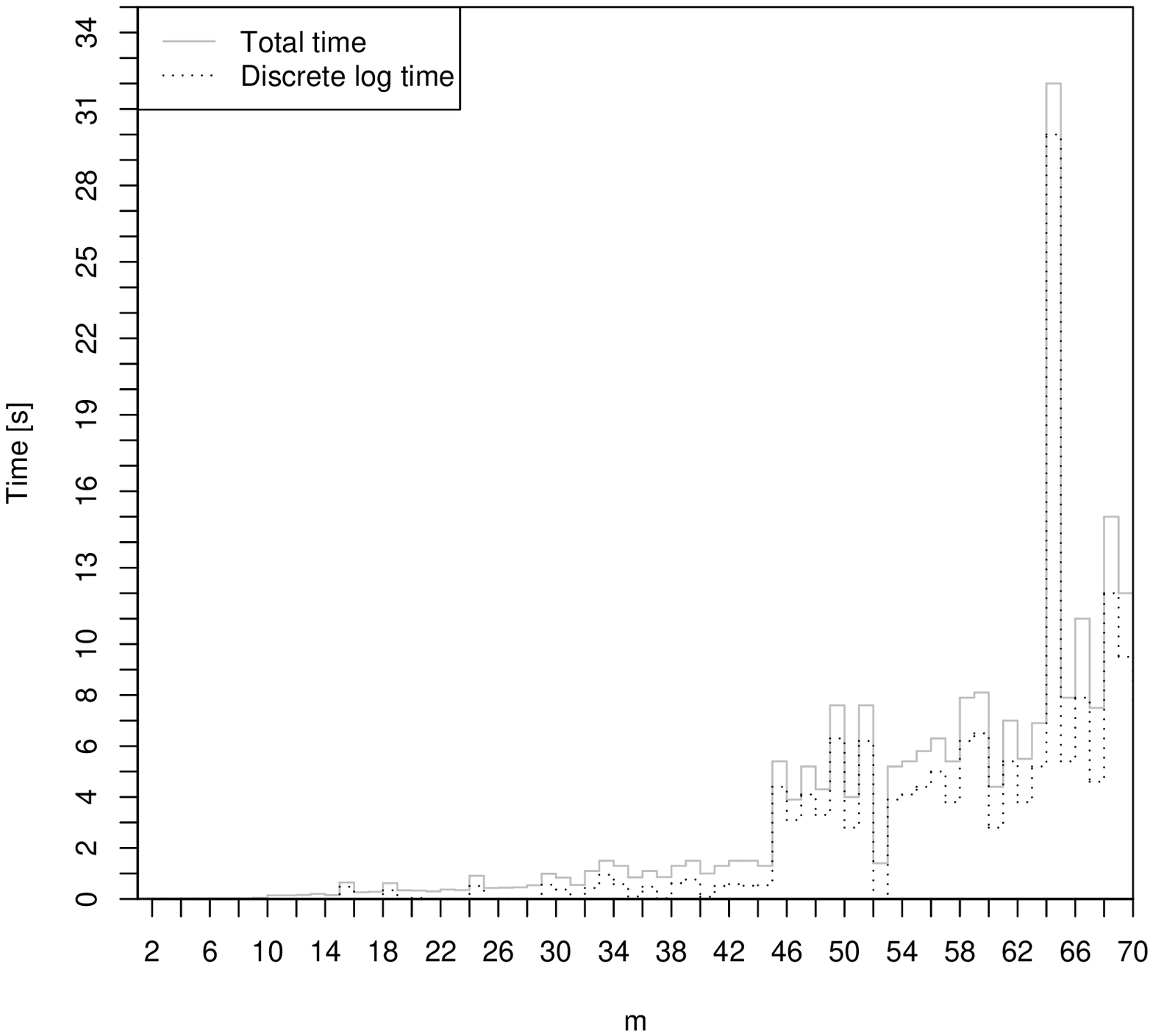}
\caption{Benchmark of stabiliser computation}
\label{fig:stab_benchmark}
\end{figure}

We have benchmarked the computation of generating sets for
stabilisers, for various field sizes, as shown in Figure
\ref{fig:stab_benchmark}. For each field size, $q = 2^{2m + 1}$,
generating sets for the stabilisers of $100$ random points were
computed, and the average running time for each call is listed. The
amount of this time that was spent in discrete logarithm computations
is also indicated.

We used the software package \textsc{R} (see \cite{r_man}), to produce
Figures \ref{fig:recognition_benchmark},
\ref{fig:conjugation_benchmark} and \ref{fig:stab_benchmark}.

All benchmarks were carried out using $\MAGMA$ V2.12-9, on a PC with
an Intel Xeon CPU running at $2.8$ GHz and with $1$ GB of RAM.  For
the conjugation problem, the highest value of $m$ was $55$, since
higher field sizes required too much memory. For the recognition and
stabiliser computation, there was never any shortage of memory, and
the benchmark indicated that much larger fields should also be
feasible. The expectation was that the conjugation problem and the
stabiliser computation would be much more time consuming than the
recognition, and in order to shorten the total time, $100$ rather than
$200$ computations were performed for each field size. The
benchmark confirmed this expectation.

Moreover, the benchmark was also used as a way to check Conjecture
\ref{conjecture_correctness}. Each stabiliser
computation involves at least $2$ calls to Algorithm
\ref{alg:stab_element}, so at least $14000$ checks of the conjecture
was made during the benchmark. The fact that it never failed 
provides strong evidence to support the conjecture.

%\appendix

%\input{tables.tex}

%\end{titlepage}

%\tableofcontents

%\listofalgorithms

%\mainmatter

%\backmatter

\bibliographystyle{amsplain}
\bibliography{suzuki}

\end{document}